\documentclass[12pt,reqno]{amsart}

\usepackage{latexsym,amsfonts,amssymb,amsmath,amscd}
\usepackage{epsfig}
\usepackage{float}
\usepackage{url}

\newcommand{\Z}{{\mathbb Z}}
\newcommand{\Q}{{\mathbb Q}}
\newcommand{\C}{{\mathbb C}}
\newcommand{\R}{{\mathbb R}}

\newcommand{\GG}{{\mathcal G}}

\newcommand{\www}{\widetilde}

\newcommand{\oooo}{\overline}

\DeclareMathOperator{\id}{id}

\DeclareMathOperator{\rank}{rank}

\DeclareMathOperator{\Rad}{Rad}

\begin{document}

\theoremstyle{plain}
\newtheorem{lemma}{Lemma}[section]
\newtheorem{definition/lemma}[lemma]{Definition/Lemma}
\newtheorem{lemma/definition}[lemma]{Lemma/Definition}
\newtheorem{theorem}[lemma]{Theorem}
\newtheorem{proposition}[lemma]{Proposition}
\newtheorem{corollary}[lemma]{Corollary}
\newtheorem{conjecture}[lemma]{Conjecture}
\newtheorem{conjectures}[lemma]{Conjectures}

\theoremstyle{definition}
\newtheorem{definition}[lemma]{Definition}
\newtheorem{withouttitle}[lemma]{}
\newtheorem{remark}[lemma]{Remark}
\newtheorem{remarks}[lemma]{Remarks}
\newtheorem{example}[lemma]{Example}
\newtheorem{examples}[lemma]{Examples}
\newtheorem{notations}[lemma]{Notations}

\title[Presentations of Weyl group elements]
{Reduced and nonreduced presentations of Weyl group elements} 

\author{Sven Balnojan \and Claus Hertling}

\address{Sven Balnojan\\
Lehrstuhl f\"ur Mathematik VI, Universit\"at Mannheim, Seminargeb\"aude
A 5, 6, 68131 Mannheim, Germany}

\email{sbalnoja@mail.uni-mannheim.de}

\address{Claus Hertling\\
Lehrstuhl f\"ur Mathematik VI, Universit\"at Mannheim, 
Seminargeb\"aude
A 5, 6, 68131 Mannheim, Germany}

\email{hertling@math.uni-mannheim.de}

\date{April 27, 2016}

\subjclass[2010]{17B22, 20F55}

\keywords{root system, subroot lattice, reduced presentation, 
quasi Coxeter element, extended affine root system}

\thanks{This work was supported by the DFG grant He2287/4-1
(SISYPH)}

\begin{abstract}
This paper is a sequel to work of Dynkin on subroot lattices 
of root lattices and to work of Carter on presentations 
of Weyl group elements as products of reflections.

The quotients $L/L_1$ are calculated for all irreducible
root lattices $L$ and all subroot lattices $L_1$.
The reduced (i.e. those with minimal number of reflections) 
presentations of Weyl group elements
as products of arbitrary reflections are classified. Also nonreduced
presentations are studied. Quasi Coxeter elements
and strict quasi Coxeter elements are defined and classified. 
An application to extended affine root lattices is given.
A side result is that any set of roots which generates
the root lattice contains a $\Z$-basis of the root lattice.
\end{abstract}

\maketitle

\tableofcontents

\setcounter{section}{0}

\section{Introduction}\label{s1}
\setcounter{equation}{0}

This paper is a sequel to work of Dynkin \cite[\S 5]{Dy}
on subroot lattices of root lattices and to work of Carter \cite{Ca}
on presentations of Weyl group elements as products of arbitrary reflections
(that means, at all possible roots, not only at roots of a fixed root basis 
as in \cite{Hu}).

A {\it root lattice} is a $\Z$-lattice $L$ together with a 
scalar product $(.,.):L_\R\times L_\R\to\R$ on the underlying
real vector space and a finite set $\Phi\subset L-\{0\}$ of {\it roots}
such that $\Phi$ is a generating set of $L$ as a $\Z$-lattice,
the reflection $s_\alpha:L_\R\to L_\R$ at a root $\alpha$ maps $\Phi$
to itself, $2(\beta,\alpha)/(\alpha,\alpha)\in\Z$ for any two
roots $\alpha$ and $\beta$, and $\Phi\cap\R\alpha=\{\pm\alpha\}$
for any root $\alpha$ (definition \ref{t2.2}).
Then the group $W:=\langle s_\alpha\, |\, \alpha\in \Phi\rangle
\subset O(L,(.,.))$ is the Weyl group of the root lattice.
The set $\Phi$ of roots is the {\it root system}.

Root lattices turn up in the theory of semisimple Lie algebras.
A standard reference is \cite[ch. VI]{Bo}. The irreducible 
root lattices are classified and form the series
$A_n,B_n,C_n,D_n$ and the exceptional cases $E_6,E_7,E_8,F_4$ and $G_2$.
Section \ref{s2} recalls their classification, standard models
and Dynkin diagrams. Because the condition 
$\Phi\cap\R\alpha=\{\pm\alpha\}$ for any root $\alpha$ is unnatural
in the context of a generalization of root lattices, we will consider
the slightly more general notion without this condition and call
it {\it p.n. root lattice}. Here {\it p.n.} stands for {\it possibly nonreduced}.
The classification of the irreducible p.n. root lattices contains
besides the irreducible root lattices only the series $BC_n$. 

Dynkin \cite[\S 5]{Dy} classified all isomorphism classes 
of pairs $(L,L_1)$ where $L$ is an irreducible root lattice and $L_1$ is a
subroot lattice (definition \ref{t3.1}). Crucial is an algorithm
which was rediscovered by Borel and de Siebenthal \cite{BS} 
and which allows to construct
by a sequence of two types of steps subroot lattices of a root lattice,
by adding roots to and dropping roots from a given set of roots.
We call the steps (BDdS1) and (BDdS2). In theorem \ref{t3.1}
we recall Dynkin's classification and give the easy extension to the
series $BC_n$. In theorem \ref{t3.8} we carry out the algorithm explicitly
for the exceptional cases
and construct thus for these cases for any isomorphism class of pairs $(L,L_1)$ a
representative. This allows to calculate also the quotient groups $L/L_1$
and to show the following. Define the numbers
\begin{eqnarray}\label{1.1}
k_1(L,L_1)&:=& \min(k\, |\, \textup{the group }L/L_1
\textup{ has }k\textup{ generators}),\\
k_2(L,L_1)&:=& \min(k\, |\, \exists\ 
\alpha_1,...,\alpha_k\in\Phi\textup{ s.t. }
L=L_1+\sum_{i=1}^k\Z\cdot\alpha_j),\hspace*{0.5cm}\label{1.2}\\
k_3(L,L_1)&:=& \min(k\, |\, L_1\textup{ can be constructed with }
k\\
&&\hspace*{1.5cm}\textup{of the steps (BDdS1) and (BDdS2)}).
\nonumber
\label{1.3}
\end{eqnarray}

\begin{theorem}\label{t1.1} (Part of theorem \ref{t3.8})
\begin{eqnarray}\label{1.4}
k_1(L,L_1)=k_2(L,L_1)=k_3(L,L_1).
\end{eqnarray}
\end{theorem}

Here the inequalities $k_1(L,L_1)\leq k_2(L,L_1)\leq k_3(L,L_1)$
are quite obvious. The inverse inequalities require the knowledge
of the groups $L/L_1$ and an explicit execution of the algorithm
with the steps (BDdS1) and (BDdS2).
The tables 3.1--3.6 contain  Dynkin's classification, the groups
$L/L_1$ and the numbers $k_1(L,L_1)$. 
For the series, the execution of the algorithm is less
important than for the exceptional cases, as for the series
one can associate graphs $\GG(A)$ to sets $A\subset\Phi$ of roots
which are helpful for understanding such sets and the subroot lattices
which they generate (definition \ref{t3.5} and lemma \ref{t3.6}).

Section \ref{s4} gives a proof of the following 
basic fact which seems to have been unnoticed
up to now and which may be of some independent interest. 

\begin{theorem}\label{t1.2} (Theorem \ref{t4.1})
Let $(L,(.,.),\Phi)$ be a p.n. root lattice.
Let $A\subset\Phi$ be any set of roots which generates the 
lattice $L$ as a $\Z$-module.
Then $A$ contains a $\Z$-basis of $L$.
\end{theorem}

The proof for the series is easy, it uses the graphs $\GG(A)$ 
(lemma \ref{t4.2}).
The proof for the exceptional cases is a case discussion.
Thanks to the results in section \ref{s3}, it can be reduced to
a discussion of just a few cases, but they require some detailed work
(lemma \ref{t4.3} and lemma \ref{t4.4}).
Theorem \ref{t4.1} is crucial for the proof of theorem \ref{t6.2}.

By definition, any Weyl group element $w$ can be written as a 
product of reflections at roots,
$$w=s_{\alpha_1}\circ ...\circ s_{\alpha_k},\quad \alpha_1,...,\alpha_k\in\Phi.$$
Then the tuple $(\alpha_1,...,\alpha_k)$ is called a {\it presentation}
of $w$, $k$ is its length, 
and the subroot lattice $L_1:=\sum_{i=1}^k\Z\alpha_i\subset L$
is called the subroot lattice of this presentation (definition \ref{t5.1}). 
The {\it length} $l(w)$ is the minimum of the lengths of all presentations of $w$.
A presentation is {\it reduced} if its length is $l(w)$.
Carter \cite[lemma 2 and lemma 3]{Ca} proved the following.

\begin{lemma}\label{t1.3} (Part of lemma \ref{t5.2})
A presentation $(\alpha_1,...,\alpha_k)$ 
of a Weyl group element $w$ is reduced if and only if $\alpha_1,...,
\alpha_k$ are linearly independent. And then the subroot lattice $L_1$
of the presentation satisfies
\begin{eqnarray}\label{1.5}
L_{1,\C}:=\bigoplus_{j=1}^{l(w)}\C\alpha_j = \bigoplus_{\lambda\neq 1}
\ker(w-\lambda\id:L_\C\to L_\C).\
\end{eqnarray}
\end{lemma}

\begin{definition}\label{t1.4} (Definition \ref{t5.3})
A Weyl group element $w$ is a {\it quasi Coxeter element} if it
has a reduced presentation whose subroot lattice is the full lattice $L$.
It is a {\it strict quasi Coxeter element} if the subroot lattice
of any reduced presentation is the full lattice $L$.
\end{definition}

In the homogeneous cases, quasi Coxeter elements and strict quasi 
Coxeter elements agree. In these cases, the definition of 
quasi Coxeter elements is due to Voigt \cite[Def. 3.2.1]{Vo}.
In the inhomogeneous cases, definition \ref{t1.4} is new.
Theorem \ref{t5.6} gives the classification of the quasi Coxeter
elements and the strict quasi Coxeter elements for all irreducible
p.n. root lattices. In the homogeneous cases this is an easy consequence
of the results of Carter \cite{Ca}. But in the inhomogeneous cases
and especially in the case $F_4$, there is some additional
work to do (lemma \ref{t5.7}).

Any Weyl group element can be written as a product
of strict quasi Coxeter elements for a suitable orthogonal sum
of irreducible subroot lattices. 
The tables 7-11 in \cite{Ca} give for any conjugacy class only
one subroot lattice $L_1$ and only one presentation as a 
strict quasi Coxeter element for this subroot lattice. 

Different presentations of one Weyl group element $w$ 
may have different subroot lattices $L_1$ and $L_1'$. 
Only $L_{1,\Q}=L_{1,\Q}'$ is clear, due to lemma \ref{t1.3}.
Theorem \ref{t5.10} complements \cite{Ca} and 
gives for the exceptional cases all subroot lattices of presentations
and all presentations as quasi Coxeter elements for these subroot lattices.
Here the case $F_4$ is more difficult than the cases $E_6,E_7$ and $E_8$.

The theorems \ref{t5.10}, \ref{t5.6} and \ref{t3.1} together allow to
recover the complete classification of all conjugacy classes
of Weyl group elements for the irreducible p.n. root lattices in \cite{Ca}
and provide additional information.

Section \ref{s6} studies nonreduced presentations of Weyl group elements.
Define the numbers
\begin{eqnarray}\label{1.6}
k_4(L,w)&:=&\min(k_2(L,L_1)\, |\, 
\textup{a reduced presentation of } w\\
&& \hspace*{2cm}\textup{with subroot lattice }L_1\textup{ exists}), \nonumber \\
k_5(L,w)&:=&\min(k\, |\, \textup{a presentation }(\alpha_1,...,\alpha_{l(w)+2k})
\textup{ with}\label{1.7}\\
&&\hspace*{2cm}\textup{subroot lattice the full lattice exists}). \nonumber
\end{eqnarray}
It is easy to see $k_5(L,w)\leq k_4(L,w)$.

\begin{theorem}\label{t1.5} (Theorem \ref{t6.2})
\begin{eqnarray}\label{1.8}
k_5(L,w)=k_4(L,w).
\end{eqnarray}
\end{theorem}

The proof builds on theorem \ref{t4.1} and some additional arguments
especially for the cases $C_n$ and $F_4$. 

Theorem \ref{t6.2} has an application to {\it extended affine root lattices}
in section \ref{s7}.
They had been defined by K. Saito \cite[(1.2) and (1.3)]{Sa}, 
see also \cite{AABGP}\cite{Az} and definition \ref{t7.1}.
One simply replaces in the definition of a p.n. root lattice the scalar
product by a positive semidefinite bilinear form
$(.,.):L_\Q\times L_\Q\to\Q$.
Then the quotient $L/\Rad(L)$ becomes in a natural way a p.n. root lattice.
Any element $w\in W(L)$ induces an element $\oooo{w}\in W(L/\Rad(L))$.
Presentations and quasi Coxeter elements in an extended affine root lattice
are defined as in a p.n. root lattice. 
The simple lemma \ref{t7.4} gives for a quasi Coxeter element $w\in W(L)$ 
the inequalities
\begin{eqnarray}\label{1.9}
l(\oooo{w})&\leq& \rank L-\rank\Rad(L),\\
l(\oooo{w})+2k_5(L/\Rad(L),\oooo{w})&\leq &\rank L.\label{1.10}
\end{eqnarray}
\eqref{1.10} gives a constraint on the elements $\oooo{w}$ which are
induced by quasi Coxeter elements. Theorem \ref{t6.2} says
$k_5(L/\Rad(L),\oooo{w})=k_4(L/\Rad(L),\oooo{w})$, and theorem \ref{t5.10}
allows to calculate this number.

\section{Basic facts on (possibly nonreduced) root lattices}
\label{s2}
\setcounter{equation}{0}

\noindent
This section recalls some basic facts on root systems.
A standard reference is \cite[ch. VI]{Bo}. Though
we follow the more recent notations and call 
{\it root systems} what is call there {\it reduced root systems}.
We call {\it p.n. root systems} ({\it p.n.} for 
{\it possibly nonreduced}) what is called there {\it root systems}.
We include the p.n. root lattices because the condition
\eqref{2.7} below, which distinguishes root systems,
is not necessarily preserved if one goes from an 
{\it extended affine root lattice} (see section \ref{s7}) 
to a quotient lattice.

\begin{notations}\label{t2.1}
(i) A free $\Z$-module $L$ of rank $n\in\Z_{>0}$ 
is called a {\it lattice}.
Then $L_\Q:=L\otimes_\Z\Q$, $L_\R:=L\otimes_\Z\R$ and 
$L_\C:=L\otimes_\Z\C$.

Let $L$ be a lattice and $(.,.)$ be scalar product on $L_\R$.
For $\alpha\in L-\{0\}$ and $\beta\in L$ define
\begin{eqnarray}\label{2.1}
\langle\beta,\alpha\rangle:=\frac{2(\beta,\alpha)}{(\alpha,\alpha)}.
\end{eqnarray}
Then 
\begin{eqnarray}\label{2.2}
s_\alpha:L_\R\to L_\R,\quad s_\alpha(x):=x-\langle x,\alpha\rangle
\cdot\alpha \end{eqnarray}
is a reflection. Two reflections $s_\alpha$ and $s_\beta$ satisfy
\begin{eqnarray}\label{2.2b}
s_\alpha s_\beta =s_\beta s_{s_\beta(\alpha)} = s_{s_\alpha(\beta)}s_\alpha.
\end{eqnarray}
\end{notations}

\begin{definition}\label{t2.2}
(a) A {\it p.n. root lattice} is a triple $(L,(.,.),\Phi)$
where $L$ is a lattice, $(.,.):L_\R\times L_\R\to\R$ is 
a scalar product, and $\Phi\subset L-\{0\}$ is a finite set 
such that the following properties hold.
\begin{eqnarray}\label{2.3}
&&\Phi\textup{ is a generating set of }L\textup{ as a }
\Z\textup{-module}.\\ \label{2.4}
&&\textup{For any }\alpha\in\Phi\ s_\alpha(\Phi)=\Phi.\\ \label{2.5}
&&\langle\beta,\alpha\rangle\in\Z\textup{ for any }
\alpha,\beta\in\Phi.
\end{eqnarray}
The elements of $\Phi$ are the {\it roots}, and $\Phi$
is a {\it p.n. root system}. The finite group
\begin{eqnarray}\label{2.6}
W:=\langle s_\alpha\, |\, \alpha\in \Phi\rangle
\subset O(L,(.,.))
\end{eqnarray}
is the {\it Weyl group}.

(b) A {\it root lattice} is a p.n. root lattice
$(L,(.,.),\Phi)$ which satisfies additionally the condition:
\begin{eqnarray}\label{2.7}
\textup{If }\alpha\in\Phi,\textup{ then }\Phi\cap\R\alpha
=\{\pm\alpha\}.
\end{eqnarray}
Then $\Phi$ is a {\it root system}.

\medskip
(c) (Lemma) The orthogonal sum of several (p.n.) root lattices
is (in a most natural way) a (p.n.) root lattice.

\medskip
(d) A (p.n.) root lattice is {\it irreducible} 
if it is not isomorphic to the orthogonal sum of 
several (p.n.) root lattices.
\end{definition}

The classification of p.n. root lattices and of root lattices
is as follows. Again, a standard reference is \cite[ch. VI]{Bo}.

\begin{theorem}\label{t2.3}
(a) Any (p.n.) root lattice is either irreducible or 
isomorphic to an orthogonal sum of several irreducible
(p.n.) root lattices. 

\medskip
(b) If $(L,(.,.),\Phi)$ is an irreducible (p.n.) root lattice
then also $(L,c\cdot(.,.),\Phi)$ for any $c\in\R_{>0}$
is an irreducible (p.n.) root lattice. Two irreducible
(p.n.) root lattices are of the same {\it type} if they
differ up to isomorphism only by such a scalar $c$.

\medskip
(c) The types of irreducible p.n. root lattices
are given by 5 series and 5 exceptional ones with the 
following names,
\begin{eqnarray}\label{2.8}
A_n\ (n\geq 1),\ B_n\ (n\geq 2),\ C_n\ (n\geq 3),\ BC_n\ (n\geq 1),\\ 
D_n\ (n\geq 4), E_6,\ E_7,\ E_8,\ F_4,\ G_2.\nonumber
\end{eqnarray}
All except $BC_n$ are root lattices.

\medskip
(d) The following list presents one irreducible p.n. root lattice 
of each type. Always $L_\R\subset\R^m$ for some $m\in\{n,n+1,n+2\}$.
Here $(.,.)$ is the restriction to $L_\R$ of the 
standard scalar product on $\R^m$, 
and $e_1,...,e_m$ is the standard ON-basis of $\R^m$.
\begin{eqnarray} \label{2.9}
{\bf A_n:} & m=n+1,& 
\Phi=\{\pm(e_i-e_j)\, |\, 1\leq i<j\leq n+1\}.\\
{\bf B_n:} & m=n,& 
\Phi=\{\pm e_i\, |\, 1\leq i\leq n\}\label{2.10}\\
&&\cup\ \{\pm e_i\pm e_j\, |\, 1\leq i<j\leq n\}.\nonumber\\
{\bf C_n:} & m=n,&
\Phi=\{\pm e_i\pm e_j\, |\, 1\leq i<j\leq n\}\label{2.11}\\
&&\cup\ \{\pm 2e_i\, |\, 1\leq i\leq n\}.\nonumber\\
{\bf BC_n:} & m=n,&
\Phi=\{\pm e_i\, |\, 1\leq i\leq n\}\label{2.12}\\
&&\cup\ \{\pm e_i\pm e_j\, |\, 1\leq i<j\leq n\}
\cup\{\pm 2e_i\, |\, 1\leq i\leq n\}.\nonumber\\
{\bf D_n:} & m=n,&
\Phi=\{\pm e_i\pm  e_j\, |\, 1\leq i<j\leq n\}.
\label{2.13}
\end{eqnarray}
\begin{eqnarray}
{\bf E_6:} & m=8,&
\Phi=\{\pm e_i\pm e_j\, |\, 3\leq i<j\leq 7\}.
\label{2.14}\\
&&\cup\ \{\frac{1}{2}\sum_{i=1}^8\varepsilon_ie_i\, |\, 
\varepsilon_i=\pm 1,\varepsilon_1=\varepsilon_2=\varepsilon_8,
\prod_{i=1}^8 \varepsilon_i=1\}\nonumber\\
{\bf E_7:} & m=8,&
\Phi=\{\pm e_i\pm e_j\, |\, 2\leq i<j\leq 7\}
\cup\ \{\pm (e_1+e_8)\}\label{2.15} \\
&&\cup\ \{\frac{1}{2}\sum_{i=1}^8\varepsilon_ie_i\, |\, 
\varepsilon_i=\pm 1,\varepsilon_1=\varepsilon_8,
\prod_{i=1}^8 \varepsilon_i=1\}.\nonumber\\ 
{\bf E_8:} & m=8,&
\Phi=\{\pm e_i\pm e_j\, |\, 1\leq i<j\leq 8\}
\label{2.16}\\
&&\cup\ \{\frac{1}{2}\sum_{i=1}^8\varepsilon_ie_i\, |\, 
\varepsilon_i=\pm 1, \prod_{i=1}^8 \varepsilon_i=1\}.\nonumber\\
{\bf F_4:} & m=4,&
\Phi=\{\pm e_i\, |\, 1\leq i\leq 4\}\label{2.17}\\
&&\cup\ \{\pm e_i\pm e_j\, |\, 1\leq i<j\leq 4\}
\cup\{\frac{1}{2}(\pm e_1\pm e_2\pm e_3\pm e_4)\}.\nonumber\\
{\bf G_2:} & m=3,&
\Phi=\{\pm (e_i-e_j)\, |\, 1\leq i<j\leq 3\}
\label{2.18}\\
&&\cup\ \{\pm(2e_{\pi(1)}-e_{\pi(2)}-e_{\pi(3)}\, |\, \pi\in S_3\}.
\nonumber
\end{eqnarray}

\end{theorem}

\begin{remarks}\label{t2.4}
(i) The p.n. root lattices above have roots of the following lengths,
\begin{eqnarray*}
\begin{array}{c|c|c|c|c|c|c|c}
A_n&D_n&E_n&B_n&F_4&C_n&G_2&BC_n\\ 
\sqrt{2}&\sqrt{2}&\sqrt{2}&1,\sqrt{2}&1,\sqrt{2}&
\sqrt{2},2&\sqrt{2},\sqrt{6}&1,\sqrt{2},2
\end{array}
\end{eqnarray*}
The root lattices of types $A_n,D_n,E_n$ have only roots
of one length and are therefore called {\it homogeneous}.
The root lattices of types $B_n,C_n,F_4$ and $G_2$ have roots
of two lengths, {\it short} and {\it long} roots.
The p.n. root lattices $BC_n$ have roots of three lengths,
short, long and {\it extra long} roots.

\medskip
(ii) In the tables 3.1 -- 3.4 the symbols $A_n,...,G_2$ will denote
root lattices with roots of lengths as above. There we will also
consider a few root systems with other lengths, and a few
other names for some of the root lattices above:
\begin{eqnarray*}
\begin{array}{l}
A_0=B_0=BC_0=\{0\}\textup{ denotes the rank 0 lattice.}\\
D_2:=2A_1:=A_1\perp A_1,\quad D_3:=A_3.\\
\www A_1=B_1\textup{ denotes a root lattice of type }A_1
\textup{ with roots of length }1.\\
C_1\textup{ denote a root lattice of type }A_1
\textup{ with roots of length }2.\\
C_2\textup{ denotes a root lattice of type }B_2
\textup{ with roots of lengths }\sqrt{2}\textup{ and }2.
\end{array}
\end{eqnarray*}
In the table 3.5 the roots in the root systems of type $C_3$ have
lengths $1$ and $\sqrt{2}$. In the table 3.6, roots in 
$A_2$ and $A_1$ have length $\sqrt{6}$, roots in $\www A_1$
have length $\sqrt{2}$.

\medskip
(iii) The Weyl group $W(A_n)$ of the root lattice above of type
$A_n$ acts on the basis
$e_1,...,e_{n+1}$ of $\R^{n+1}\supset L_\R$ by permutations,
$W(A_n)\cong S_{n+1}$, and $\sigma\in S_{n+1}$ maps
$e_i$ to $e_{\sigma(i)}$.

The Weyl groups of the p.n. root lattices above of the types 
$B_n,C_n$ and $BC_n$ coincide and act on the basis
$e_1,...,e_n$ of $\R^n=L_\R$ by {\it signed permutations},
$$W(B_n)=W(C_n)=W(BC_n)\cong \{\pm 1\}^n\rtimes S_n,$$
and $(\varepsilon_1,...,\varepsilon_n,\sigma)\in 
\{\pm 1\}^n\rtimes S_n$ maps $e_i$ to $\varepsilon_i e_{\sigma(i)}$.

The Weyl group of the root lattice above of type $D_n$
is the subgroup of index 2 given by the condition
$\prod_{i=1}^n\varepsilon_i=1$.

\medskip
(iv) Let $(L,(.,.),\Phi)$ be an irreducible root lattice.
To any subset $A=\{\delta_1,...,\delta_l\}\subset\Phi$
with $A\cap (-A)=\emptyset$
we associate a {\it generalized Dynkin diagram} as follows.
It is a graph with $l$ vertices, labelled $\delta_1,...,\delta_l$.
Between vertices $\delta_i$ and $\delta_j$ with $i\neq j$
there is no edge or an edge with additional information as follows.
\begin{eqnarray*}
\textup{no edge} && \textup{if }(\delta_i,\delta_j)=0\\
\textup{a normal edge} && \textup{if }
\|\delta_i\|=\|\delta_j\|\textup{ and }
\langle\delta_i,\delta_j\rangle =-1,\\
\textup{a dotted edge} && \textup{if }
\|\delta_i\|=\|\delta_j\|\textup{ and }
\langle\delta_i,\delta_j\rangle=1,\\
\textup{a double arrow from }\delta_i\textup{ to }\delta_j && \textup{if }
\|\delta_i\|=\sqrt{2}\|\delta_j\|\textup{ and }
\langle\delta_i,\delta_j\rangle=-2,\\
\textup{a double dotted arrow from }\delta_i\textup{ to }\delta_j && \textup{if }
\|\delta_i\|=\sqrt{2}\|\delta_j\|\textup{ and }
\langle\delta_i,\delta_j\rangle=2,\\
\textup{a triple arrow from }\delta_i\textup{ to }\delta_j && \textup{if }
\|\delta_i\|=\sqrt{3}\|\delta_j\|\textup{ and }
\langle\delta_i,\delta_j\rangle=-3
\end{eqnarray*}
The corresponding pictures are depicted below.
\begin{figure}[!h]
	\centering
\includegraphics[scale=1]{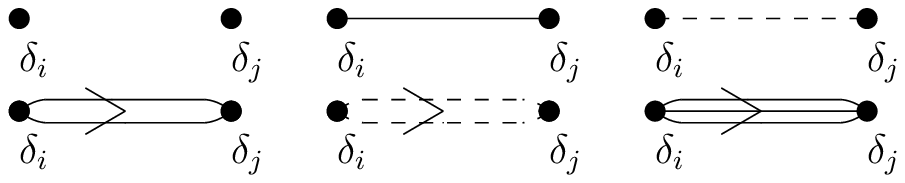}
\end{figure}

\medskip
Other cases will not be considered.
If $A$ is a $\Z$-basis of $L$ and the diagram is connected, 
then the diagram encodes
up to a common scalar the intersection numbers $(\delta_i,\delta_j)$,
and thus it determines the irreducible root system.

\medskip
(v) The following list gives for each of the root lattices in theorem
\ref{t2.3} (d) a {\it root basis} $\delta_1,...,\delta_n$ 
(a $\Z$-basis of $L$ with additional properties \cite{Bo})
and an additional root $\delta_{n+1}$ (which is minus
the maximal root with respect to the root basis).
The diagram for the root basis is called {\it Dynkin diagram},
the diagram for $\delta_1,...,\delta_{n+1}$ is called 
{\it extended Dynkin diagram}. 
The roots $\delta_1,...,\delta_{n+1}$ satisfy a linear
relation. For the cases $E_6,E_7,E_8,F_4,G_2$, it is given in
lemma \ref{t3.7} (c). In the case of $E_6$, 
$\delta_7=\frac{1}{2}(-\sum_{i=1,2,3,8}e_i
+\sum_{i=4,5,6,7}e_i).$
\begin{eqnarray*}
\begin{array}{l|l|l}
\textup{type} & \delta_1,...,\delta_n & \delta_{n+1}\\ \hline
A_n & e_i-e_{i+1}\ (i=1,...,n) & -e_1+e_{n+1} \\
B_n & -e_1,\ e_i-e_{i+1}\ (i=1,..,n-1) & e_{n-1}+e_n\\
C_n & -2e_1,\ e_i-e_{i+1}\ (i=1,..,n-1) & 2e_n\\
D_n & e_i-e_{i+1}\ (i=1,...,n-1),\ e_{n-1}+e_n & -e_1-e_2\\
E_6 & \frac{1}{2}\sum_{i=1}^8e_i,\ -e_3-e_4,\ 
e_i-e_{i+1}\ (i=3,...,6) & \delta_7\\
E_7 & \frac{1}{2}\sum_{i=1}^8e_i,\ -e_2-e_3,\ 
e_i-e_{i+1}\ (i=2,...,6) & -e_1-e_8\\ 
E_8 & \frac{1}{2}\sum_{i=1}^8e_i,\ -e_1-e_2,\ 
e_i-e_{i+1}\ (i=1,...,6) & e_7-e_8\\
F_4 & \frac{1}{2}\sum_{i=1}^4e_i,\ -e_1,\ 
e_1-e_2,\ e_2-e_3 & e_3-e_4\\
G_2 & e_1-e_2 ,\ -e_1+2e_2-e_3 & -e_1-e_2+2e_3
\end{array}
\end{eqnarray*}
The following table gives the extended Dynkin diagrams.

\begin{figure}[!h]
	\centering
	\includegraphics[scale=1]{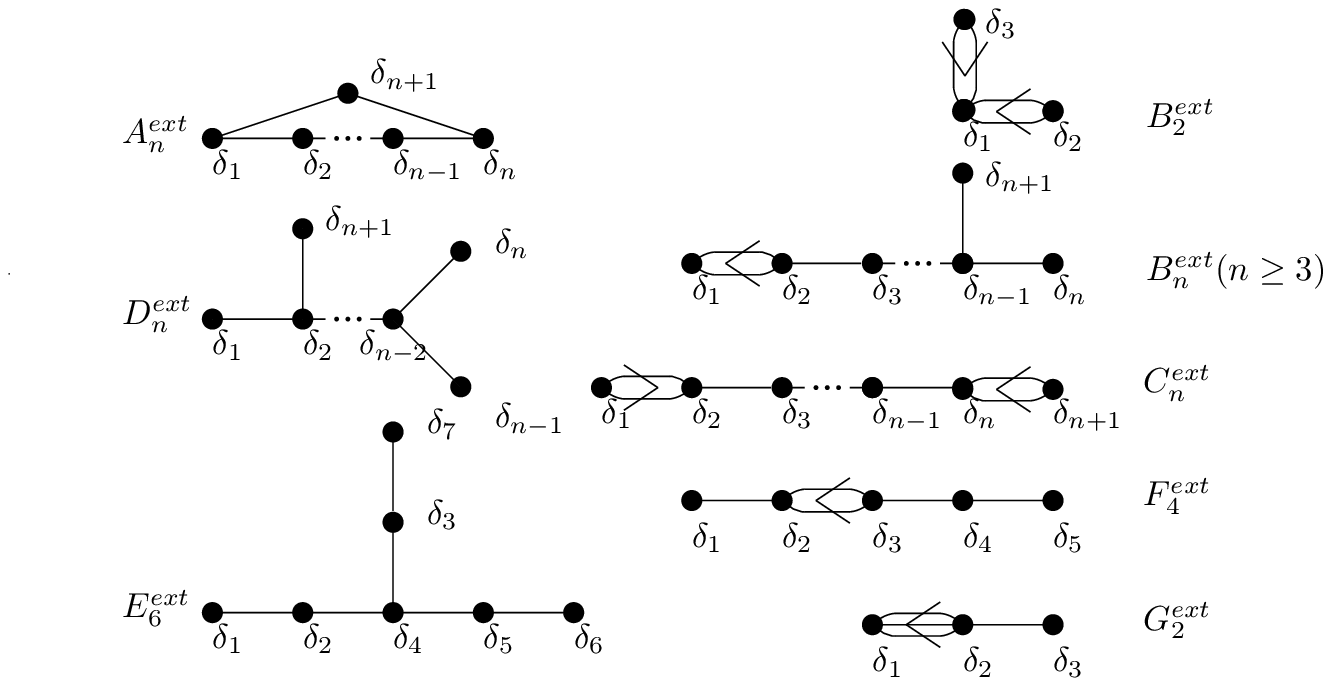}
\end{figure}
\begin{figure}[!h]
	\centering
	\includegraphics[scale=1]{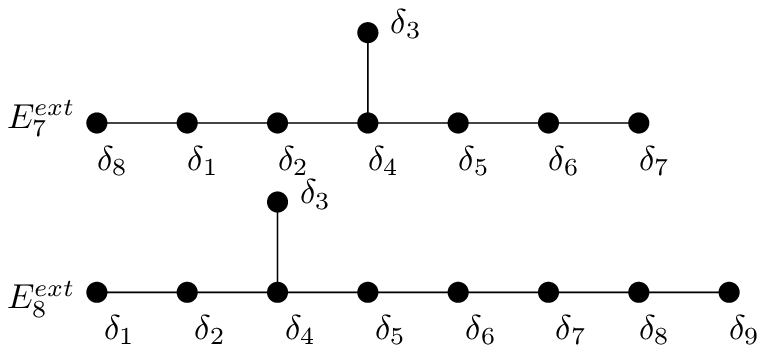}
\end{figure}

\noindent
The Dynkin diagram encodes
up to a common scalar the intersection numbers $(\delta_i,\delta_j)$
of the basis $\delta_1,...,\delta_n$, 
and thus it determines the irreducible root system.
\end{remarks}

\section{Subroot lattices and quotients}\label{s3}
\setcounter{equation}{0}

\noindent
The subroot lattices of an irreducible root lattice
can be determined up to isomorphism by a recipe due to 
Dynkin \cite{Dy} and Borel and de Siebenthal \cite{BS}.
In \cite{Dy} also a list and additional information is given.
In this section,
we will review the recipe and go through it, and thus we will 
recover the list. We will extend the list in two ways.
First, we consider also the p.n. subroot lattices of the
p.n. root lattices of type $BC_n$. Second, we will calculate
for any (isomorphism class of a) pair $(L,L_1)$ where $L_1$ is a 
(p.n.) subroot lattice of an irreducible (p.n.) root lattice $L$
the quotient group $L/L_1$. This will be helpful in section
\ref{s4} and crucial in the sections 
\ref{s6} and \ref{s7}.

\begin{definition}\label{t3.1}
Let $(L,(.,.),\Phi)$ be a (p.n.) root lattice.

\medskip
(a) A (p.n.) root lattice $(L_1,(.,.)_1,\Phi_1)$ is a 
{\it (p.n.) subroot lattice} of $(L,(.,.),\Phi)$ if
$L_1\subset L$ and $(.,.)_1$ is the restriction of $(.,.)$
to $L_1$ and $\Phi_1=L_1\cap \Phi$.

A notation: Because $(.,.)_1$ and $\Phi_1$ are determined
by $L_1$, we will talk of the subroot lattice $L_1$.

\medskip
(b) A (p.n.) root lattice $(L_1,(.,.)_1,\Phi_1)$ is the
(p.n.) root lattice of a {\it (p.n.) subroot system} if $L_1\subset L$ 
and $(.,.)_1$ is the restriction of $(.,.)$ to $L_1$ and 
$\Phi_1\subset L_1\cap \Phi$.

\medskip
(c) The index of a subroot lattice $L_1$ is $[L\cap L_{1,\Q}:L_1]\in\Z_{\geq 1}$.
\end{definition}

\begin{remarks}\label{t3.2}
Let $(L,(.,.),\Phi)$ be a (p.n.) root lattice.

\medskip
(i) Let $L_1\subset L$ be a $\Z$-sublattice. 
Define $(.,.)_1$ as the restriction of $(.,.)$ to $L_1$.
Define $\Phi_1:=L_1\cap\Phi$. Then
$(L_1,(.,.)_1,\Phi_1)$ is a (p.n.) subroot lattice
if and only if it is a (p.n.) root lattice, and this holds
if and only if $L_1$ is generated by $\Phi_1$ as a $\Z$-module:
\eqref{2.5} holds for $\Phi_1$, and $s_\alpha(\Phi_1)\subset L_1$
holds for $\alpha\in\Phi_1$ because of \eqref{2.2} and \eqref{2.1}.
This gives \eqref{2.4}.

\medskip
(ii) If $A=\{\alpha_1,...,\alpha_l\}\subset\Phi$ is any nonempty subset, then the data 
\begin{eqnarray}\label{3.1}
L_1:=\sum_{i=1}^l\Z\cdot\alpha_i,\quad
(.,.)_1:=(.,.)_{|L_1},\quad \Phi_1:=L_1\cap\Phi
\end{eqnarray}
satisfy the conditions in (i) and are a (p.n.) subroot lattice.
%

\medskip
(iii) Any (p.n.) subroot lattice is the root lattice of a 
subroot system. If $(L,(.,.),\Phi)$ is a homogeneous root lattice
also the inverse holds. But if $(L,(.,.),\Phi)$ contains orthogonal
summands which are of types $B_k,C_k,BC_k,F_4$ or $G_2$, then there
are subroot lattices $(L_1,(.,.),\Phi_1)$ such that the subsets
$\Phi_2\subsetneqq\Phi_1$ of short roots give rise to root lattices
$(L_1,(.,.),\Phi_2)$ of subroot systems $\Phi_2$, such that these
root lattices are not subroot lattices.
We will not work much with them, but in \cite{Ca} they are used.

\medskip
(iv) If one erases from any of the extended Dynkin diagrams
one vertex, one obtains a disjoint union of Dynkin diagrams.
This leads to the following recipe with two kinds of steps 
with which one obtains
easily subroot lattices of a root lattice. It is due to
\cite[\S 5]{Dy} and \cite{BS}, therefore we call the steps
(BDdS1) and (BDdS2). Start with a root lattice $(L,(.,.),\Phi)$.
Choose a root basis $A\subset\Phi$, 
that is a $\Z$-basis of $L$ consisting of roots 
such that its generalized Dynkin diagram 
(defined in remark \ref{t2.4} (iv)) is a disjoint union of Dynkin diagrams.
$L$ decomposes uniquely into an orthogonal sum of irreducible
subroot lattices, which are called the {\it summands} of $L$.

\begin{list}{}{}
\item[Step (BDdS1):] Choose one summand
$L_1$ of $L$, add to $A$ the unique root $\www\delta$ in $\Phi_1$ which
gives together with the roots in $A\cap\Phi_1$ an extended
Dynkin diagram (it is a linear combination of the roots
in $A\cap\Phi_1)$ and delete from $A\cup\{\www\delta\}$
an arbitrary root in $A\cap\Phi_1$. The new set $\www A\subset\Phi$
defines a subroot lattice $\www L$ of $L$ of the same rank as $L$.
\item[Step (BDdS2):] Choose one summand
$L_1$ of $L$ and delete from $A$
an arbitrary root in $A\cap\Phi_1$. The new set $\www A\subset\Phi$
defines a subroot lattice $\www L$ of $L$ with $\rank L_1=\rank L-1$.
\end{list}

In both cases $\www A$ is a root basis of $\www L$.
Therefore one can repeat the steps.
The change in the Dynkin diagrams is easy to see.
In the step (BDdS1) one extends one component to its extended version
and then erases one vertex. In the step (BDdS2) one simply erases
one vertex.
\end{remarks}

The following theorem is mainly due to Dynkin 
\cite[\S 5]{Dy},
the recipe in part (a) is also in \cite{BS}.
The only new (though rather trivial) part is the discussion of 
the cases $BC_n$. That will follow from lemma \ref{t3.6} 
below.

\begin{theorem}\label{t3.3}
(a) Let $(L,(.,.),\Phi)$ be a root lattice.
Any subroot lattice is obtained by the choice of a suitable
root basis of $L$ and by a suitable sequence of the
steps (BDdS1) and (BDdS2).

\medskip
(b) The first columns of the tables 3.1 -- 3.6 
list all isomorphism classes of
pairs $((L,(.,.),\Phi),L_1)$ where $(L,(.,.),\Phi)$ is
an irreducible (p.n.) root lattice with the lengths of
the roots as in theorem \ref{t2.3} (d) and where $L_1$ is
a subroot lattice. 

The tables give the name for the type of $L_1$,
where additionally the lengths of the roots of the summands of 
$L_1$ are taken into account. The symbols
$A_0,B_0,BC_0,D_2,D_3,\www A_1,B_1,C_1,C_2$ from remark
\ref{t2.4} (ii) are used.
The new notations $[...]'$ and $[...]''$ are explained in (d) below.

\medskip
(c) With one class of exceptions, the following holds.
If $((L,(.,.),\Phi),L_1)$ and $((L,(.,.),\Phi),L_2)$ are
isomorphic pairs as in (b), then a Weyl group element $w\in W$
with $w(L_1)=w(L_2)$ exists.
The class of exceptions are the sublattices of $D_n$ 
of types $A_{k_1}+...+A_{k_r}$ with all $k_1,...,k_r$ odd.
For each of those types there are two conjugacy classes
with respect to $W$.

\medskip
(d) The tables 3.3 and 3.4 contain pairs $[H]'$ and $[H]''$ with
$H\in\{A_5+A_1,A_5,A_3+2A_1,A_3+A_1,4A_1,3A_1\}$ for $E_7$
and with $H\in\{A_7,A_5+A_1,2A_3,A_3+2A_1,4A_1\}$ for $E_8$.
Here $[H]'$ and $[H]''$ denote (classes in the sense of (b) of)
subroot lattices which are isomorphic
if one forgets the embedding into $L$. 
But for a subroot lattice $L_1\subset L$ of type $[H]'$ and a 
subroot lattice $L_2\subset L$ of type $[H]''$, 
the pairs $(L,L_1)$ and $(L,L_2)$ are not isomorphic.
This is an implication of the following properties:
A subroot lattice $L_3\subset L$ of type $A_7$ for $E_7$ and
of type $A_8$ for $E_8$ with $L_1\subset L_3\subset L$ exists,
but no subroot lattice $L_4\subset L$ of type $A_7$ for $E_7$ and
of type $A_8$ for $E_8$ with $L_2\subset L_4\subset L$ exists.
\end{theorem}

The informations in the following tables 3.1 -- 3.6
are treated in theorem \ref{t3.3}, lemma \ref{t3.7} 
and theorem \ref{t3.8}. Always $L$ is one of the 
p.n. subroot lattices in theorem \ref{t2.3} (d),
and $L_1$ is a p.n. subroot lattice of the type
indicated. In the tables 3.2--3.6 it is
$L_1:=\sum_{i\in\{1,...,n\}\cup I-J}\Z\cdot\delta_i$.
Here the roots $\delta_k$ for $k\geq n+2$ (in the cases $E_7,E_8,F_4$)
are defined in Lemma \ref{t3.7}.
The quotient $L/L_1$ is given up to isomorphism.
Here $\Z_m:=\Z/m\Z$ for $m\in\Z_{>0}$.
For $k_1=k_1(L,L_1)$ see theorem \ref{t3.8}.
For the symbols $A_0,B_0,BC_0,D_2,D_3,\www A_1,B_1,C_1,C_2$ 
see remark \ref{t2.4} (ii).

\medskip
{\bf Table 3.1 for $A_n,B_n,C_n,BC_n,D_n$:} Here $r\geq 0$, $s\geq 0$,
$a_i\geq 0$, $b_j\geq 1$ in the cases $C_{b_j}$, 
$b_j\geq 2$ in the cases $D_{b_j}$,
$m=-1$ in the case $A_n$,
$m=0$ in the cases $C_n$ and $D_n$,
$m\geq 0$ in the cases $B_n$ and $BC_n$,
\begin{eqnarray}
\label{3.2}
\sum_{i=1}^r (a_i+1)+\sum_{j=1}^s b_j +m&=&n.
\end{eqnarray}
\begin{eqnarray*}
\begin{array}{l|l|l|ll}
L & L_1 & L/L_1 & k_1(L,L_1) \\ \hline
A_n & \sum_{i=1}^r A_{a_i} & \Z^{r-1} & r-1&\\
B_n & \sum_{i=1}^r A_{a_i}+\sum_{j=1}^s D_{b_j} + B_m & 
\Z^r\times \Z^s_2 & r+s&\\
C_n & \sum_{i=1}^r A_{a_i}+\sum_{j=1}^s C_{b_j} & 
\Z^r\times \Z^{s-1}_2 & r+s-1&\textup{if }s\geq 1\\
 & & \Z^r & r&\textup{if }s=0\\
BC_n & \sum_{i=1}^r A_{a_i}+\sum_{j=1}^s C_{b_j} +BC_m& 
\Z^r\times \Z^s_2 & r+s&\\
D_n & \sum_{i=1}^r A_{a_i}+\sum_{j=1}^s D_{b_j} & 
\Z^r\times\Z^{s-1}_2 & r+s-1 &\textup{if }s\geq 1\\
 & & \Z^r & r &\textup{if }s=0
\end{array}
\end{eqnarray*}

\medskip
{\bf Table 3.2 for $E_6$:} 
\begin{eqnarray*}
\begin{array}{l|l|l|l|l}
L_1 & I & J & L/L_1 & k_1\\ \hline
E_6 & - & - & \{0\} & 0 \\
A_5+A_1 & 7 & 2 & \Z_2 & 1 \\
3A_2 & 7 & 4 & \Z_3 & 1 \\ \hline
A_5 & - & 3 & \Z & 1 \\
2A_2+A_1 & - & 4 & \Z & 1 \\
A_4 + A_1 & - & 2 & \Z & 1 \\
D_5 & - & 1 & \Z & 1 \\ 
A_3 +2A_1 & 7 & 2,3 & \Z\times \Z_2 & 2 \\ \hline
A_4 & - & 1,2 & \Z^2 & 2 \\
A_3+A_1 & - & 2,3 & \Z^2 & 2
\end{array}
\hspace*{0.5cm}
\begin{array}{l|l|l|l|l}
L_1 & I & J & L/L_1 & k_1\\ \hline
2A_2 & - & 3,4 & \Z^2 & 2 \\
A_2+2A_1 & - & 4,5 & \Z^2 & 2 \\
4A_1 & 7 & 2,3,5 & \Z^2\times \Z_2 & 3 \\
D_4 & - & 1,6 & \Z^2 & 2 \\ \hline
A_3 & - & 1,2,3 & \Z^3 & 3 \\
A_2+A_1 & - & 2,3,6 & \Z^3 & 3 \\
3A_1 & - & 2,3,5 & \Z^3 & 3 \\ \hline
A_2 & - & 1,2,3,4 & \Z^4 & 4 \\
2A_1 & - & 1,2,3,5 & \Z^4 & 4 \\ \hline
A_1 & - & 1,...,5 & \Z^5 & 5
\end{array}
\end{eqnarray*}

\medskip
{\bf Table 3.3 for $E_7$:} 
\begin{eqnarray*}
\begin{array}{l|l|l|l|l}
L_1 & I & J & L/L_1 & k_1 \\ \hline
E_7 & - & - & \{0\} & 0 \\
D_6+A_1 & 8 & 6 & \Z_2 & 1 \\
A_5+A_2 & 8 & 5 & \Z_3 & 1 \\
2A_3+A_1 & 8 & 4 & \Z_4 & 1 \\
A_7 & 8 & 3 & \Z_2 & 1 \\
D_4+3A_1 & 8,9 & 1,6 & \Z^2_2 & 2 \\
7A_1 & 8,9,10 & 1,6,4 & \Z^3_2 & 3 \\ \hline
E_6 & - & 7 & \Z & 1 \\
D_5+A_1 & - & 6 & \Z & 1 \\
A_4+A_2 & - & 5 & \Z & 1 
\end{array}
\hspace*{0.5cm}
\begin{array}{l|l|l|l|l}
L_1 & I & J & L/L_1 & k_1 \\ \hline
A_3\! +\! A_2\! +\! A_1 & - & 4 & \Z & 1 \\
{}[A_5+A_1]' & 8 & 1,3 & \Z\times\Z_2 & 2 \\
{}[A_5+A_1]'' & - & 2 & \Z & 1 \\
D_6 & - & 1 & \Z & 1 \\
D_4+2A_1 & 8 & 1,6 & \Z\times\Z_2 & 2 \\
A_3+3A_1 & 8 & 1,4 & \Z\times\Z_2 & 2 \\
3A_2 & 8 & 2,5 & \Z\times\Z_3 & 2 \\
2A_3 & 8 & 3,4 & \Z\times\Z_2 & 2 \\
A_6 & - & 3 & \Z & 1 \\
6A_1 & 8,9 & 1,4,6 & \Z\times\Z^2_2 & 3 \\ \hline
\end{array}
\end{eqnarray*}
\begin{eqnarray*}
\begin{array}{l|l|l|l|l}
D_5 & - & 1,7 & \Z^2 & 2 \\
A_4+A_1 & - & 3,6 & \Z^2 & 2 \\
2A_2+A_1 & - & 2,5 & \Z^2 & 2 \\
{[}A_5]' & - & 1,3 & \Z^2 & 2 \\
{[}A_5]'' & - & 1,2 & \Z^2 & 2 \\
D_4+A_1 & - & 1,6 & \Z^2 & 2 \\
A_3+A_2 & - & 3,5 & \Z^2 & 2 \\
5A_1 & 8 & 1,4,6 & \Z^2\times\Z_2 & 3 \\
A_2+3A_1 & - & 4,6 & \Z^2 & 2 \\
{[}A_3+2A_1]' & 8 & 1,3,4 & \Z^2\times\Z_2 & 3 \\
{[}A_3+2A_1]'' & - & 1,4 & \Z^2 & 2 \\ \hline
D_4 & - & 1,6,7 & \Z^3 & 3 \\
A_4 & - & 1,2,3 & \Z^3 & 3 
\end{array}
\hspace*{0.5cm}
\begin{array}{l|l|l|l|l}
{[}A_3+A_1]' & - & 1,3,4 & \Z^3 & 3 \\
{[}A_3+A_1]'' & - & 1,2,4 & \Z^3 & 3 \\
2A_2 & - & 1,2,5 & \Z^3 & 3 \\
A_2+2A_1 & - & 2,4,5 & \Z^3 & 3 \\
{[}4A_1]' & 8 & 1,3,4,6 & \Z^3\times\Z_2 & 4 \\
{[}4A_1]'' & - & 2,4,6 & \Z^3 & 3 \\ \hline
A_3 & - & 1,2,3,4 & \Z^4 & 4 \\
A_2+A_1 & - & 1,2,3,5 & \Z^4 & 4 \\
{[}3A_1]' & - & 1,3,4,6 & \Z^4 & 4 \\
{[}3A_1]'' & - & 1,2,4,6 & \Z^4 & 4 \\ \hline
A_2 & - & 1,..,5 & \Z^5 & 5 \\
2A_1 & - & 1,..,4,6 & \Z^5 & 5 \\ \hline
A_1 & - & 1,..,6 & \Z^6 & 6 
\end{array}
\end{eqnarray*}

\medskip
{\bf Table 3.4 for $E_8$:} 
\begin{eqnarray*}
\begin{array}{l|l|l|l|l}
L_1 & I & J & L/L_1 & k_1 \\ \hline
E_8 & - & - & \{0\} & 0 \\
A_8 & 9 & 3 & \Z_3 & 1 \\
D_8 & 9 & 1 & \Z_2 & 1 \\
A_7+A_1 & 9 & 2 & \Z_4 & 1 \\
A_5\! +\! A_2\! +\! A_1 & 9 & 4 & \Z_6 & 1 \\
2A_4 & 9 & 5 & \Z_5 & 1 \\
4A_2 & 9,13 & 7,4 & \Z_3^2 & 2 
\end{array}
\hspace*{0.5cm}
\begin{array}{l|l|l|l|l}
L_1 & I & J & L/L_1 & k_1 \\ \hline
E_6+A_2 & 9 & 7 & \Z_3 & 1 \\
E_7+A_1 & 9 & 8 & \Z_2 & 1 \\
D_6+2A_1 & 9,10 & 1,8 & \Z_2^2 & 2 \\
D_5+A_3 & 9 & 6 & \Z_4 & 1\\
2D_4 & 9,10 & 1,6 & \Z_2^2 & 2 \\
D_4+4A_1 & 9,10,12 & 1,6,4 & \Z_2^3 & 3 \\
2A_3\! +\! 2A_1 & 9,10 & 8,4 & \Z_2\times \Z_4 & 2 \\
8A_1 & 9,10,11,12 & 1,6,8,4 & \Z_2^4 & 4 \\ \hline
\end{array}
\end{eqnarray*}
\begin{eqnarray*}
\begin{array}{l|l|l|l|l}
A_6+A_1 & - & 2 & \Z & 1\\
A_4\! +\! A_2\! +\! A_1 & - & 4 & \Z & 1\\
A_5+A_2 & 9 & 3,4 & \Z\times\Z_3 & 2\\
3A_2+A_1 & 9 & 4,7 & \Z\times\Z_3 & 2\\
E_6+A_1 & - & 7 & \Z & 1\\
E_7 & - & 8 & \Z & 1\\
D_7 & - & 1 & \Z & 1\\
D_5+2A_1 & 9 & 6,8 & \Z\times\Z_2 & 2\\
D_4+3A_1 & 9,10 & 1,4,6 & \Z\times\Z_2^2 & 3\\
2A_3+A_1 & 9 & 2,6 & \Z\times\Z_4 & 2
\end{array}
\hspace*{0.5cm}
\begin{array}{l|l|l|l|l}
7A_1 & 9,10,11 & 1,6,8,4 & \Z\times\Z_2^3 & 4\\
D_6+A_1 & 9 & 1,8 & \Z\times\Z_2 & 2\\
D_5+A_2 & - & 6 & \Z & 1\\
A_3\! +\! A_2\! +\! 2A_1 & 9 & 4,6 & \Z\times\Z_2 & 2\\
D_4+A_3 & 9 & 1,6 & \Z\times\Z_2 & 2\\
A_3+4A_1 & 9,10 & 1,4,8 & \Z\times\Z_2^2 & 3\\
A_4+A_3 & - & 5 & \Z & 1\\
A_5+2A_1 & 9 & 1,4 & \Z\times\Z_2 & 2\\
{}[A_7]' & - & 3 & \Z & 1\\
{}[A_7]'' & 9 & 1,2 & \Z\times\Z_2 & 2 \\  \hline 
\end{array}
\end{eqnarray*}
\begin{eqnarray*}
\begin{array}{l|l|l|l|l}
3A_2 & 9 & 3,4,7 & \Z^2\times\Z_3 & 3\\
E_6 & - & 7,8 & \Z^2 & 2\\
D_6 & - & 1,8 & \Z^2 & 2\\
D_4+2A_1 & 9 & 1,6,8 & \Z^2\times\Z_2 & 3\\
{}[2A_3]' & - & 3,5 & \Z^2 & 2\\
{}[2A_3]'' & 9 & 1,2,6 & \Z^2\times\Z_2 & 3\\
D_5+A_1 & - & 1,7 & \Z^2 & 2\\
A_3+3A_1 & 9 & 1,4,6 & \Z^2\times\Z_2 & 3\\
D_4+A_2 & - & 1,6 & \Z^2 & 2
\end{array}
\hspace*{0.5cm}
\begin{array}{l|l|l|l|l}
6A_1 & 9,10 & 1,4,6,8 & \Z^2\times\Z_2^2 & 4\\
A_2+4A_1 & 9 & 4,6,8 & \Z^2\times\Z_2 & 3\\
A_4+2A_1 & - & 1,4 & \Z^2 & 2\\
A_6 & - & 1,3 & \Z^2 & 2\\
A_3\! +\! A_2\! +\! A_1 & - & 4,8 & \Z^2 & 2\\
{}[A_5+A_1]' & - & 2,3 & \Z^2 & 2\\
{}[A_5+A_1]'' & 9 & 1,2,4 & \Z^2\times\Z_2 & 3\\
A_4+A_2 & - & 3,4 & \Z^2 & 2\\
2A_2+2A_1 & - & 4,6 & \Z^2 & 2 \\ \hline
\end{array}
\end{eqnarray*}
\begin{eqnarray*}
\begin{array}{l|l|l|l|l}
D_5 & - & 1,7,8 & \Z^3 & 3\\
{}[A_3+2A_1]' & - & 2,3,5 & \Z^3 & 3\\
{}[A_3+2A_1]'' & 9 & 1,2,4,6 & \Z^3\times\Z_2 & 4\\
A_3+A_2 & - & 3,5,8 & \Z^3 & 3\\
A_5 & - & 1,2,8 & \Z^3 & 3\\
5A_1 & 9 & 1,4,6,8 & \Z^3\times\Z_2 & 4\\
A_4+A_1 & - & 1,3,4 & \Z^3 & 3\\
D_4+A_1 & - & 1,6,7 & \Z^3 & 3\\
A_2+3A_1 & - & 1,4,6 & \Z^3 & 3\\
2A_2+A_1 & - & 2,3,6 & \Z^3 & 3 \\ \hline
\end{array}
\hspace*{0.5cm}
\begin{array}{l|l|l|l|l}
D_4 & - & 1,6,7,8 & \Z^4 & 4\\
{}[4A_1]' & - & 2,3,5,7 & \Z^4 & 4\\
{}[4A_1]'' & 9 & 1,2,4,6,8 & \Z^4\times\Z_2 & 5\\
A_2\! +\! 2A_1 & - & 1,2,4,6 & \Z^4 & 4\\
2A_2 & - & 1,2,3,6 & \Z^4 & 4\\
A_3+A_1 & - & 1,2,3,5 & \Z^4 & 4\\
A_4 & - & 1,2,3,4 & \Z^4 & 4\\   \hline
A_3 & - & 1,2,3,4,5 & \Z^5 & 5\\
A_2+A_1 & - & 1,2,3,4,6 & \Z^5 & 5\\
3A_1 & - & 1,2,3,5,7 & \Z^5 & 5 \\ \hline
\end{array}
\end{eqnarray*}
\begin{eqnarray*}
\begin{array}{l|l|l|l|l}
A_2 & - & 1,..,6 & \Z^6 & 6\\
2A_1 & - & 1,..,5,7 & \Z^6 & 6\\ \hline
\end{array}
\hspace*{0.5cm}
\begin{array}{l|l|l|l|l}
A_1 & - & 1,..,7 & \Z^7 & 7 
\end{array}
\end{eqnarray*}

\medskip
{\bf Table 3.5 for $F_4$:} 
\begin{eqnarray*}
\begin{array}{l|l|l|l|l}
L_1 & I & J & L/L_1 & k_1 \\ \hline
F_4 & - & - & \{0\} & 0 \\
B_4 & 5 & 1 & \Z_2 & 1 \\
A_3+\www A_1 & 5 & 2 & \Z_4 & 1 \\
A_2+\www A_2 & 5 & 3 & \Z_3 & 1 \\
C_3+A_1 & 5 & 4 & \Z_2 & 1 \\
D_4 & 5,6 & 1,2 & \Z_2^2 & 2 \\
B_2+2A_1 & 5,6 & 1,4 & \Z_2^2 & 2 \\
4A_1 & 5,6,7 & 1,2,4 & \Z_2^3 & 3 \\ \hline
B_3 & - & 1 & \Z & 1 \\
B_2+A_1 & 5 & 1,4 & \Z\times\Z_2 & 2 \\
A_2+\www A_1 & - & 2 & \Z & 1 
\end{array}
\hspace*{0.5cm}
\begin{array}{l|l|l|l|l}
A_3 & 5 & 1,2 & \Z\times\Z_2 & 2 \\
2A_1+\www A_1 & 5 & 2,4 & \Z\times\Z_2 & 2 \\
A_1+\www A_2 & - & 3 & \Z & 1 \\
C_3 & - & 4 & \Z & 1 \\
3A_1 & 5,6 & 1,2,4 & \Z\times\Z_2^2 & 3 \\ \hline
A_2 & - & 1,2 & \Z^2 &  2\\
B_2 & - & 1,4 & \Z^2 &  2\\
A_1+\www A_1 & - & 2,3 & \Z^2 & 2 \\
2A_1 & 5 & 1,2,4 & \Z^2\times\Z_2 & 3 \\
\www A_2 & - & 3,4 & \Z^2 & 2 \\ \hline
\www A_1 & - & 2,3,4 & \Z^3 & 3 \\
A_1 & - & 1,2,3 & \Z^3 & 3 
\end{array}
\end{eqnarray*}

\medskip
{\bf Table 3.6 for $G_2$:} 
\begin{eqnarray*}
\begin{array}{l|l|l|l|l}
L_1 & I & J & L/L_1 & k_1 \\ \hline
G_2 & - & - & \{0\} & 0 \\
A_2 & 3 & 1 & \Z_3& 1 
\end{array}
\hspace*{0.5cm}
\begin{array}{l|l|l|l|l}
A_1+\www A_1 & 3 & 2 & \Z_2 & 1\\ \hline
\www A_1 & - & 2 & \Z & 1 \\
A_1 & - & 1 & \Z & 1
\end{array}
\end{eqnarray*}

\begin{remarks}\label{t3.4}
(i) Let $(L,(.,.),\Phi)$ be an irreducible root lattice.
Theorem \ref{t3.3} (a)+(b)+(d) tells the following.
There is an almost 1:1 correspondence between 
the set of isomorphism classes of pairs $((L,(.,.),\Phi),L_1)$ 
with $L_1$ a subroot lattice and 
the set of unions of Dynkin diagrams which are obtained by 
iterations of the
graphical versions of the steps (BDdS1) and (BDdS2) in remark
\ref{t3.2} (iv), namely 

(BDdS1): Go from one Dynkin diagram to the extended Dynkin diagram
and erase an arbitrary vertex.

(BDdS2): Erase an arbitrary vertex.

The only exceptions are the pairs $[H]'$ and $[H]''$ discussed
in theorem \ref{t3.3} (d). They have the same Dynkin diagrams.

\medskip
(ii) In the table 11 in \cite[ch. II, \S 5]{Dy} there are
two misprints. $A_6+A_2$ has to be replaced by $E_6+A_2$.
And one of the two $A_7+A_1$ has to be replaced by $E_7+A_1$.
\end{remarks}

In the cases of the series $A_n,B_n,C_n$ and $D_n$,
one can see the subroot lattices also in a different way,
by associating a graph to a generating set $A\subset\Phi_1$
of a subroot lattice $L_1$. This works also in the case
of the series $BC_n$ and will give the proof of the statements
in theorem \ref{t3.3} for $BC_n$. The graphs are defined as follows.

\begin{definition}\label{t3.5}
Let $(L,(.,.),\Phi)$ be a p.n. root lattice 
in theorem \ref{t2.3} (d) of one of the
types $A_n,B_n,C_n,BC_n,D_n$.
Let $A=\{\alpha_1,...,\alpha_l\}\subset\Phi$ be a nonempty subset.
It defines a p.n. subroot lattice $L_1=\sum_{i=1}^l\Z\cdot \alpha_i$.
A graph $\GG(A)$ with or without markings of the vertices 
and with one or two types of edges is defined as follows.

\medskip
(a) $L$ of type $A_n$: The graph $\GG(A)$ has $n+1$ vertices
which are labelled $1,...,n+1$. It has $l$ edges. A root
$\alpha\in A$ with $\alpha=\pm(e_i-e_j)$ gives an edge between
the vertices $i$ and $j$. So, if $e_i-e_j$ and $e_j-e_i$
are in $A$, there are two edges between the vertices $i$ and $j$.
The same applies in the cases (b) and (c).

\medskip
(b) $L$ of type $BC_n$: The graph $\GG(A)$ has $n$ vertices
which are labelled $1,...,n$. Any root $\pm e_i$ in $A$
leads to a marking of the vertex $i$ which is called a 
{\it short} marking (and which may be represented by
a circle around the vertex). 
Any root $\pm 2e_i$ in $A$ leads to a marking of the vertex $i$
which is called a {\it long} marking (and which may be represented
by a square around the vertex). So, depending on how many of
the roots $\pm e_i$ and $\pm 2 e_i$ are in $A$, the vertex
$i$ has between 0 and 4 markings.
Any root $\pm(e_i-e_j)$ gives a normal edge between the
vertices $i$ and $j$. Any root $\pm (e_i+e_j)$ gives a 
dotted edge between the vertices $i$ and $j$. 
So, between the vertices $i$ and $j$ there are between 0 and 4
edges.

\medskip
(c) $L$ of type $B_n$ or $C_n$ or $D_n$: The graph is
defined as in the case of type $BC_n$.
(In the case of $B_n$ there are no long markings,
in the case of $C_n$ there are no short markings,
in the case of $D_n$ there are no markings at all).
\end{definition}

The following lemma is obvious.

\begin{lemma}\label{t3.6}
Consider the same data as in definition \ref{t3.5}.
The orthogonal irreducible summands of the subroot lattice $L_1$
can be read off from the graph $\GG(A)$ as follows.
Each of the following subgraphs yields a summand,
which is generated by the roots which contribute via
markings or edges to this subgraph.

\begin{list}{}{}
\item[${\bf A_k}$:] 
A component of $\GG(A)$ which has no markings
and in which any cycle has an even number of dotted edges
yields a summand of type $A_k$. Here $k+1$ is the number
of vertices of the component. 
(An isolated vertex with no markings yields thus the 
summand $A_0=\{0\}$).
\item[${\bf B_k}$ or ${\bf BC_k}$:]
The union of all components of $\GG(A)$ which contain a vertex with
a short marking yields a summand of type $B_k$ if
$L$ is of type $B_n$ and a summand of type $BC_k$ if $L$
is of type $BC_n$. Here $k$ is the number of vertices
of the union of these components. If this union is empty,
we write $B_0(=\{0\})$ if $L$ is of type $B_n$ and 
$BC_0(=\{0\})$ if $L$ is of type $BC_n$.
\item[${\bf C_k}$:]
A component of $\GG(A)$ which does not contain a vertex with a short
marking, but which contains a cycle with an odd number
of dotted edges or which contains a vertex with a long
marking yields a summand of type $C_k$ if $L$ is of type 
$C_n$ or $BC_n$. Here $k$ is the number of vertices of this
component. 
\item[${\bf D_k}$:]
A component of $\GG(A)$ 
which does not contain a vertex with a marking,
but which contains a cycle with an odd number of dotted
edges yields a summand of type $D_k$ if $L$ is of type
$B_n$ or $D_n$. Here $k$ is the number of vertices of this
component.
\end{list}
\end{lemma}

The statements in theorem \ref{t3.3} for the cases
$A_n,B_n,C_n,BC_n$ and $D_n$ follow easily
from this lemma and from the structure of the
Weyl group, which was described in remark \ref{t2.4} (iii).

\medskip
For the calculation of the quotients $L/L_1$ in theorem
\ref{t3.8}, we need in the cases $E_6,E_7,E_8,F_4,G_2$
a concrete subroot lattice $L_1$ for each isomorphism
class of pairs $(L,L_1)$. This is found in lemma \ref{t3.7}
by carrying out the recipe with the steps
(BDdS1) and (BDdS2).

\begin{lemma}\label{t3.7}
Let $L$ be an irreducible root lattice
in theorem \ref{t2.3} (d) of one of the types
$E_6,E_7,E_8,F_4,G_2$. Additionally to the roots
$\delta_1,...,\delta_{n+1}$ which are defined in 
remark \ref{t2.4} (v), the following roots are considered.
\begin{eqnarray}\label{3.3}
\begin{array}{l|l}
E_7 & \delta_9=e_6+e_7,\ 
\delta_{10}=e_4+e_5.\\ 
E_8 & \delta_{10}=e_7+e_8,\ \delta_{11}=e_5+e_6,
\ \delta_{12}=e_3+e_4,\\
 & \delta_{13}=
\frac{1}{2}(-e_1+\sum_{i=2}^5e_i-\sum_{i=6}^8 e_i).\\
F_4 & \delta_6=e_3+e_4,\ \delta_7=-e_1-e_2.
\end{array}
\end{eqnarray}

\medskip
(a) In the cases $E_7,E_8$ and $F_4$, 
the generalized Dynkin diagrams which take into account
the roots $\delta_1,...,\delta_{n+1}$ and the roots above
look as follows.

\noindent
\begin{figure}[!h]
	\centering
	\includegraphics[scale=1]{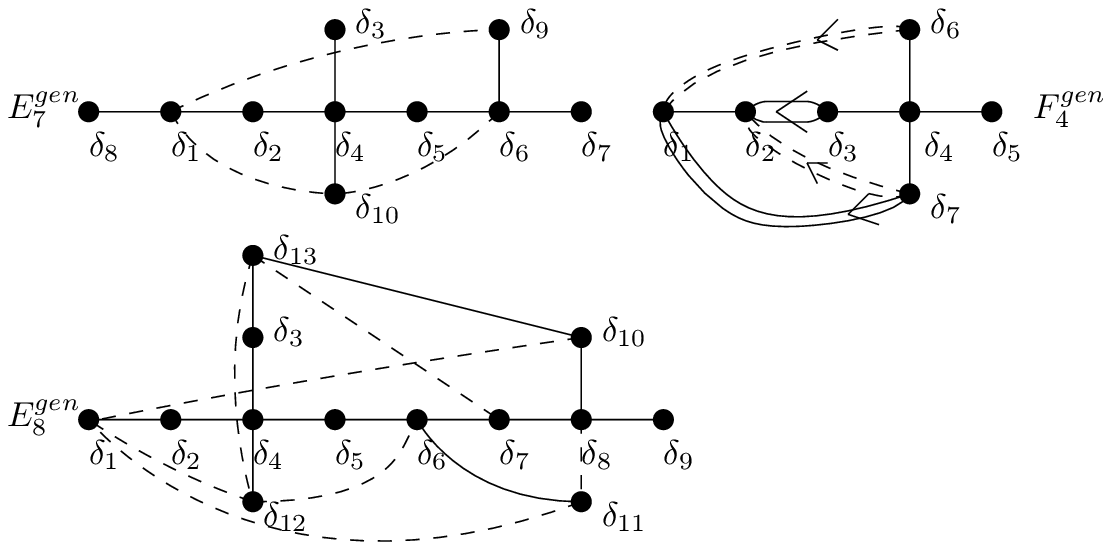}
\end{figure}

(The edges which are not horizontal or vertical will be 
irrelevant except for the dotted edge between $\delta_1$
and $\delta_{10}$ in the Dynkin diagram of $E_8$.
It will be used once, in the construction of a subroot lattice
of type $2A_3+2A_1$.)

\medskip
(b) The second and third column in the tables 3.2 -- 3.6 encode a
realization of the recipe in remark \ref{t3.2} (iv),
in the following way.

Suppose in the tables 3.2 -- 3.6 
in the line for one subroot lattice, 
$I$ is given as the sequence $i_1,...,i_r$ of numbers and 
$J$ is given as the sequence $j_1,...,j_s$ of numbers, 
with $0\leq r\leq s$. 
One carries out $r$ steps (BDdS1): In the $k$-th step one
adds the root $\delta_{i_k}$ and then erases the root
$\delta_{j_k}$.
Afterwards one carries out $s-r$ steps (BDdS2):
One erases the roots $\delta_{j_{s-r+1}},...,\delta_{j_s}$.
This leads to a subroot lattice of the type indicated
in the first column.

\medskip
(c) The roots $\delta_k$ for $k\geq n+1$ are linear combinations
of the roots $\delta_1,...,\delta_n$. The linear relations are 
as follows.
\begin{eqnarray*}
{\bf E_6:}&&0=\delta_1+2\delta_2+2\delta_3+3\delta_4+2\delta_5
+\delta_6+\delta_7.\\
{\bf E_7:}&&0=2\delta_1+3\delta_2+2\delta_3+4\delta_4+3\delta_5
+2\delta_6+\delta_7+\delta_8,\\
&&0=\delta_2+\delta_3+2\delta_4+2\delta_5
+2\delta_6+\delta_7+\delta_9,\\
&&0=\delta_2+\delta_3+2\delta_4+\delta_5
+\delta_{10}.\\
{\bf E_8:}&&0=2\delta_1+4\delta_2+3\delta_3+6\delta_4+5\delta_5
+4\delta_6+3\delta_7+2\delta_8+\delta_9,\\
&&0=2\delta_1+3\delta_2+2\delta_3+4\delta_4
+3\delta_5+2\delta_6+\delta_7-\delta_{10},\\
&&0=\delta_2+\delta_3+2\delta_4+2\delta_5
+2\delta_6+\delta_7+\delta_{11},\\
&&0=\delta_2+\delta_3+2\delta_4+\delta_5+\delta_{12},\\
&&0=\delta_1+2\delta_2+2\delta_3+3\delta_4
+2\delta_5+\delta_6+\delta_{13}.\\
{\bf F_4:}&&0=2\delta_1+4\delta_2+3\delta_3+2\delta_4+\delta_5,\\
&&0=2\delta_1+2\delta_2+\delta_3-\delta_6,\\
&&0=2\delta_2+\delta_3-\delta_7.\\
{\bf G_2:}&&0=3\delta_1+2\delta_2+\delta_3.
\end{eqnarray*}

\end{lemma}

{\bf Proof:}
The proof of the lemma is tedious as there are many cases,
but the parts (a) and (c) and most of part (b) are completely elementary. In part (b) one has to check not only
that the result has the correct Dynkin diagram, but
also that the steps (BDdS1) work, i.e. that 
one has after adding a root an extended Dynkin diagram
and that then a root of this extended Dynkin diagram is erased.
The details are left to the reader.

The only nontrivial part concerns the subroot lattices
of types $[H]'$ and $[H]''$ in the cases $E_7$ and $E_8$.
One sees that in the cases of $E_k$, $k=7,8$, 
the constructed subroot lattices of types $[H]'$ are contained in the
subroot lattice $\bigoplus_{i\in\{1,2,4,..,k+1\}}\Z\delta_i$
of type $A_k$.
The constructed subroot lattices of type $[H]''$
contain in the case $E_7$ the roots $\delta_3,\delta_5,\delta_7$
and in the case $E_8$ the roots $\delta_3,\delta_5,\delta_7,\delta_9$.
The following claim shows that the constructed subroot lattices
$[H]''$ are not contained in subroot lattices of type $A_7$
respectively $A_8$.

\medskip
{\bf Claim:}
(i) 
Let $L$ be the lattice of type $E_7$ in theorem \ref{t2.3} (d). 
There is no subroot lattice $L_1$ of type $A_7$
with $\delta_3,\delta_5,\delta_7\in L_1$.

(ii) 
Let $L$ be the lattice of type $E_8$ in theorem \ref{t2.3} (d). 
There is no subroot lattice $L_1$ of type $A_8$ 
with $\delta_3,\delta_5,\delta_7,\delta_9\in L_1$.

\medskip
{\bf Proof of the claim:}
(i) Suppose that $L_1$ is a subroot lattice of type $A_7$
with $\delta_3,\delta_5,\delta_7\in L_1$. 
These three roots generate a subroot lattice $L_2\subset L_1$
of type $3A_1$.
By theorem \ref{t3.3}, up to isomorphism there is only one pair
of type $(A_7,3A_1)$. Therefore a root $\alpha\in L_1$ with
$(\alpha,\delta_3)=-1,(\alpha,\delta_5)=0,(\alpha,\delta_7)=0$
exists. But now observe $\delta_3=e_2-e_3,\delta_5=e_4-e_5,
\delta_7=e_6-e_7$. Neither the roots in $\Phi(L)$ of type 
$\pm e_i\pm e_j$ nor the roots in $\Phi(L)$ of type 
$\frac{1}{2}\sum_{i=1}^8\varepsilon_ie_i$ can serve as a 
root $\alpha$. Contradiction.

(ii) Analogously to (i). \hfill ($\Box$)

This finishes the proof of lemma \ref{t3.7}. \hfill$\Box$

\begin{theorem}\label{t3.8}
(a) Let $(L,(.,.),\Phi)$ be an irreducible p.n. root lattice,
and let $L_1$ be a p.n. subroot lattice. 
The third column in table 3.1 and the fourth column in the 
tables 3.2 -- 3.6 gives the
isomorphism class of the quotient group $L/L_1$.

\medskip
(b) Let $(L,(.,.),\Phi)$ be a p.n. root lattice,
and let $L_1$ be a p.n. subroot lattice. 
Define the numbers
\begin{eqnarray}\label{3.4}
k_1(L,L_1)&:=& \min(k\, |\, \textup{the group }L/L_1
\textup{ has }k\textup{ generators}\},\\
k_2(L,L_1)&:=& \min(k\, |\, \exists\ 
\alpha_1,...,\alpha_k\in\Phi\textup{ s.t. }
L=L_1+\sum_{i=1}^k\Z\cdot\alpha_j\},\hspace*{0.5cm}\label{3.5}\\
k_3(L,L_1)&:=& \min(k\, |\, L_1\textup{ can be constructed with }
k\\
&&\hspace*{1.5cm}\textup{of the steps (BDdS1) and (BDdS2)}\}.
\nonumber
\label{3.6}
\end{eqnarray}
Then
\begin{eqnarray}\label{3.7}
k_1(L,L_1)=k_2(L,L_1)=k_3(L,L_1).
\end{eqnarray}
The numbers are additive, i.e. if $L=L_2+L_3$ and
$L_1=L_4+L_5$ and $L_2\supset L_4, L_3\supset L_5$ then
\begin{eqnarray}\label{3.8}
k_1(L,L_1)=k_1(L_2,L_4)+k_1(L_3,L_5).
\end{eqnarray}
The last column of the tables 3.1 -- 3.6 gives the numbers
$k_1(L,L_1)$ for the pairs with $L$ irreducible.
Minimal sequences of the steps (BDdS1) and (BDdS2)
for the cases with $L$ of type $E_6,E_7,E_8,F_4,G_2$
are given in the second and third column of the tables 3-2 -- 3.6
(see lemma \ref{t3.7}).
\end{theorem}

{\bf Proof:}
(a) First we treat the cases $A_n,B_n,C_n,BC_n$ and $D_n$.
Let $A$ and $\GG(A)$ be as in definition \ref{t3.5} and lemma \ref{t3.6},
and let $L_1\subset L$ be the corresponding subroot lattice.
Let $\GG(A)=\bigcup_{k\in K}\GG_k$ (with $1\notin K$) 
be the decomposition into
subgraphs $\GG_k$ with the properties in lemma \ref{t3.6},
let $L_k$ be the subroot lattice which corresponds to the subgraph $\GG_k$, 
and let 
$$V_k:=\bigoplus_{i\textup{ is a vertex in }\GG_k}\Z\cdot e_i.$$
Then 
$$\bigoplus_{k\in K} V_k=\bigoplus_{i=1}^m\Z\cdot e_i\supset L
\supset L_1=\bigoplus_{k\in K}L_k$$
(with $m=n+1$ for $A_n$ and $m=n$ else) and 
\begin{eqnarray*}
L/L_1\subset (\bigoplus_{i=1}^m\Z\cdot e_i)/L_1 &\cong&
\bigoplus_{k\in K}V_k/L_k.
\end{eqnarray*}
The following table lists in the second and fourth line 
the isomorphism classes of the quotients in the first column.
\begin{eqnarray*}
\begin{array}{l|l|l|l|l|l}
L &A_n & B_n & C_n & BC_n & D_n\\ \hline
\bigoplus_{i=1}^m\Z e_i/L & 
\Z & \{0\} & \Z_2 & \{0\} & \Z_2 \\ \hline \hline 
L_k &A_l & B_l & C_l & BC_l & D_l\\ \hline
V_k/L_k& \Z & \{0\} & \Z_2 & \{0\} & \Z_2 
\end{array}
\end{eqnarray*}
Finally, in the cases $C_n$ and $D_n$, for any $k\in K$
$$L\not\subset L_k\oplus\bigoplus_{j\in K-\{k\}}V_k.$$
Therefore $L/L_1$ has the isomorphism type claimed in the table 3.1.

Now we treat the cases $E_6,E_7,E_8,F_4,G_2$.
Let $L_1\subset L$ be one of the subroot lattices constructed
in lemma \ref{t3.7} using the data in the tables 3.2 -- 3.6.
Let $\delta_{i_1},...,\delta_{i_r}$ respectively
$\delta_{j_1},...,\delta_{j_s}$ with $0\leq r\leq s$
be the roots in one line in the second respectively third
column of these tables. Then
\begin{eqnarray*}
L/L_1 =\frac{(\bigoplus_{k=1}^s\Z\cdot \delta_{j_k})+L_1}{L_1}
\cong \frac{\bigoplus_{k=1}^s\Z\cdot \delta_{j_k}}
{(\bigoplus_{k=1}^s\Z\cdot \delta_{j_k})\cap L_1}.
\end{eqnarray*}
The denominator of the right hand side is a $\Z$-lattice
of rank $r$ (because $\rank L_1=n-s+r$) 
and is generated by parts of those relations
in lemma \ref{t3.7} (c) which express the roots
$\delta_{i_1},...,\delta_{i_r}$ as linear combinations of
the roots $\delta_1,...,\delta_n$. 

We give one example: $(L,L_1)$ of type $(E_8,D_4+3A_1)$,
then $(\delta_{i_1},...,\delta_{i_r})=(\delta_9,\delta_{10})$
and $(\delta_{j_1},...,\delta_{j_s})=(\delta_1,\delta_4,\delta_6)$.
The relation for $\delta_9$ gives the element
$2\delta_1+6\delta_4+4\delta_6$ of 
$(\bigoplus_{k=1}^s\Z\cdot \delta_{j_k})\cap L_1$, 
the relation for $\delta_{10}$ gives the element
$2\delta_1+4\delta_4+2\delta_6$.
Therefore here
$$\frac{L}{L_1}\cong 
\frac{\Z\cdot\delta_1\oplus\Z\cdot\delta_4\oplus\Z\cdot\delta_6}
{\Z\cdot(2\delta_1+6\delta_4+4\delta_6)\oplus
\Z\cdot(2\delta_1+4\delta_4+2\delta_6)}\cong \Z\times\Z_2^2.$$
The calculations for all other cases $(L,L_1)$ are analogous.
They are tedious as there are many cases, but elementary.

(b)
The additivity of the numbers $k_1(L,L_1),k_2(L,L_1),k_3(L,L_1)$ 
is obvious. Therefore it is sufficient to prove \eqref{3.7}
for irreducible $L$.
The last column of the tables 3.1 -- 3.6 
can be read off from the second to last column immediately.
The first of the inequalities
\begin{eqnarray}\label{3.9}
k_1(L,L_1)\leq k_2(L,L_1)\leq k_3(L,L_1)
\end{eqnarray}
is obvious.
The second inequality follows
simply from the fact that in each step (BDdS1) or (BDdS2),
one root is erased.
In the cases $A_n,B_n,C_n,BC_n$ and $D_n$, 
one easily constructs the p.n. subroot lattices in
$k_1(L,L_1)$ steps.
Therefore then $k_3(L,L_1)\leq k_1(L,L_1)$, and 
equalities hold in \eqref{3.9}.

In the cases $E_6,E_7,E_8,F_4,G_2$
one observes $|J|=k_1(L,L_1)$ in the tables 3.2 -- 3.6.
In lemma \ref{t3.7} $|J|$ steps of type (BDsD1) and (BDsD2)
are used. Therefore $k_3(L,L_1)\leq k_1(L,L_1)$, and
equalities hold in \eqref{3.9}.
\hfill$\Box$

\section{Any generating set of roots contains a $\Z$-basis}\label{s4}
\setcounter{equation}{0}

\noindent
The purpose of this section is to prove the following theorem.
It is crucial in the proof of theorem \ref{t6.2}.
But it may be also of some independent interest.
The proof in the cases $A_n,B_n,C_n,BC_n,D_n$ is an almost trivial
application of the graphs in definition \ref{t3.5} and
lemma \ref{t3.6}. The proof in the cases $E_6,E_7,E_8,F_4,G_2$ 
is more involved.

\begin{theorem}\label{t4.1}
Let $(L,(.,.),\Phi)$ be a p.n. root lattice.
Let $A\subset\Phi$ be any set of roots which generates the 
lattice $L$ as a $\Z$-module.
Then $A$ contains a $\Z$-basis of $L$.
\end{theorem}

In the case of vector spaces instead of $\Z$-modules,
the analogous statement is trivial.
For $\Z$-modules it is not true in general, that any generating
set contains a basis. For example
the lattice $\Z^2$ with standard basis $(1,0),(0,1)$ has the set
$\{(1,0),(1,2),(0,3)\}$ as generating set, but any two of 
these elements generate a proper sublattice.

The rest of this section is devoted to the proof of theorem \ref{t4.1}.
It is obviously sufficient to prove it in the cases where $L$
is an irreducible p.n. root lattice.

In the cases $A_n,B_n,C_n,BC_n,D_n$, theorem \ref{t4.1}
is an immediate consequence of the following lemma.
The lemma follows directly from lemma \ref{t3.6}.

\begin{lemma}\label{t4.2}
Let $(L,(.,.),\Phi)$ be a p.n. root lattice in theorem \ref{t2.3}
(d) of one of the types $A_n,B_n,C_n,BC_n,D_n$.
Let $L=\{\alpha_1,...,\alpha_l\}\subset\Phi$ be a nonempty subset.
The properties whether $A$ is a generating set of $L$ or a 
$\Z$-basis of $L$, will be characterized by properties of
the graph $\GG(A)$ from definition \ref{t3.5}.

(a) $L$ of type $A_n$: 
\begin{list}{}{}
\item[(i)]
$A$ generates $L$ as a $\Z$-module $\iff$
$\GG(A)$ is connected.
\item[(ii)]
$A$ is a $\Z$-basis of $L$ $\iff$ 
$\GG(A)$ is a tree.
\end{list}

(b) $L$ of type $B_n$: 
\begin{list}{}{}
\item[(i)]
$A$ generates $L$ as a $\Z$-module $\iff$
each component of $\GG(A)$ contains at least one vertex
with a (automatically short) marking.
\item[(ii)]
$A$ is a $\Z$-basis of $L$ $\iff$ 
each component of $\GG(A)$ is a tree and contains exactly
one marked vertex, and the vertex has only one
(automatically short) marking.
\end{list}

(c) $L$ of type $C_n$: 
\begin{list}{}{}
\item[(i)]
$A$ generates $L$ as a $\Z$-module $\iff$
$\GG(A)$ is connected, and it contains a marking
(automatically long) or a cycle with an odd number of
dotted edges.
\item[(ii)]
$A$ is a $\Z$-basis of $L$ $\iff$ 
either $\GG(A)$ is a tree and contains exactly one marked vertex 
and the marking is simple, or $\GG(A)$ is connected and 
contains no marking, but it contains
exactly one cycle and the cycle has an odd number of 
dotted lines.
\end{list}

(d) $L$ of type $BC_n$: 
\begin{list}{}{}
\item[(i)]
$A$ generates $L$ as a $\Z$-module $\iff$
each component of $\GG(A)$ contains at least one vertex
with a short marking.
\item[(ii)]
$A$ is a $\Z$-basis of $L$ $\iff$ 
each component of $\GG(A)$ is a tree and contains exactly
one marked vertex, and there is only one
marking, and the marking is short.
\end{list}

(e) $L$ of type $D_n$: 
\begin{list}{}{}
\item[(i)]
$A$ generates $L$ as a $\Z$-module $\iff$
$\GG(A)$ is connected, and it contains a cycle with an odd
number of dotted edges.
\item[(ii)]
$A$ is a $\Z$-basis of $L$ $\iff$ 
$\GG(A)$ is connected and contains exactly one cycle,
and the cycle has an odd number of dotted lines.
\end{list}
\end{lemma}

It rests to prove theorem \ref{t4.1} in the cases 
$G_2, F_4,E_6,E_7,E_8$.
Let $L$ be the root lattice of one of these types 
in theorem \ref{t2.3} (d).
In each of these cases it is sufficient to consider 
a generating set $A=\{\alpha_1,...,\alpha_{n+1}\}\subset\Phi$ 
with $n+1$ elements (where $n$ is the rank of the root lattice). 
The cases of bigger generating sets 
can be reduced to the case of such a set 
by an easy inductive argument.

There is an up to the sign unique linear combination
\begin{eqnarray}\label{4.1}
0=\sum_{i=1}^{n+1}\lambda_i\alpha_i\quad
\textup{with }\lambda_i\in\Z, \ 
\gcd(\lambda_1,...,\lambda_{n+1})=1,
\textup{ not all }\lambda_i=0.
\end{eqnarray}
It has to be shown that an index $j$ with $\lambda_j=\pm 1$ exists.

Denote by $L_i\subset L$, $i\in\{1,...,n+1\}$, the subroot lattice
generated by $A-\{\alpha_i\}$.
Then
\begin{eqnarray}\label{4.2}
\rank L_i<n \iff \lambda_i=0.
\end{eqnarray}
If this holds for some $i$ then by induction 
on the rank of the lattice one can conclude that $A-\{\alpha_i\}$
contains a $\Z$-basis of this subroot lattice.
Then this $\Z$-basis together with $\alpha_i$ forms a 
$\Z$-basis of $L$.

Thus suppose that all $\lambda_i\notin\{0,\pm 1\}$. Then
$$[L:L_i]=|\lambda_i|,\quad L_i+\Z\alpha_i=L,\quad
k_2(L,L_i)=1.$$
A priori, there are 20 possible cases in the tables 3.2 -- 3.6,
\begin{eqnarray*} 
\begin{array}{ccc}
\begin{array}{l|l|l}
L & L_i &  |\lambda_i| \\ \hline 
F_4 & B_4 & 2\\
 & C_3+A_1 & 2\\ 
 & A_2+\www A_2 & 3\\
 & A_3+\www A_1 & 4\\ \hline
E_6 & A_5+ A_1 & 2\\ 
 & 3A_2 & 3 \\ \hline 
E_7 & D_6+ A_1 & 2\\
 & A_7 & 2 \\
 & A_5+ A_2 & 3\\
 & 2A_3+ A_1 & 4\\
\end{array}&
\begin{array}{l|l|l}
L & L_i & |\lambda_i| \\ \hline 
G_2 & A_1+\www A_1 & 2 \\
 & A_2 & 3 \\ \hline
E_8 & D_8 & 2\\
 & E_7+A_1 & 2\\
 & E_6+A_2 & 3\\
 & A_8 & 3\\
 & D_5+A_3 & 4\\
 & A_7+A_1 & 4\\
 & 2A_4 & 5\\
 & A_5+A_2+ A_1 & 6\\
\end{array}
\end{array}
\end{eqnarray*}
In the cases $G_2,F_4,E_6,E_7$, the only possible values
for $|\lambda_i|$ are in $\{2,3,4\}$. 
The condition $\gcd(\lambda_1,...,\lambda_{n+1})=1$ tells
that at least one $j\in\{1,...,n+1\}$ with $|\lambda_j|=3$ exists.

The following more complicated argument gives the same conclusion
in the case $E_8$. Assume in the case $E_8$ that all
$|\lambda_i|$ are in $\{2,4,5,6\}$.
Define the decomposition
$$I_1:=\{i\, |\, 2|\lambda_i\},\quad
I_2:=\{i\, |\, \lambda_i=\pm 5\}$$
of $\{1,...,n+1\}$ into two disjoint subset.
Because of $\gcd(\lambda_1,...,\lambda_{n+1})=1$, both are nonempty.
Define the subroot lattices
$$\www L_1:=\sum_{i\in I_1}\Z\alpha_i,\qquad
\www L_2:=\sum_{i\in I_2}\Z\alpha_i.
$$
Then $c:=\gcd(\lambda_i\, |\, i\in I_1)\in\{2,4,6\}$, and
$c^{-1}\sum_{i\in I_1}\lambda_i\alpha_i$ is a primitive vector
in $\www L_1$, and $5^{-1}\sum_{i\in I_2}\lambda_i\alpha_i$
is a primitive vector in $\www L_2$. But
$$\sum_{i\in I_1}\lambda_i\alpha_i 
= -\sum_{i\in I_2}\lambda_i\alpha_i.$$
Therefore the order of the torsion part of $L/\www L_1$
is divisible by 5. But table 3.4 contains only one type of subroot
lattices with this property, the type $2A_4$. Therefore 
$|I_1|=8, |I_2|=1, I_2=\{j_0\}$ for some index $j_0$,
and $\frac{5}{2}\alpha_{j_0}\in \www L_1\subset L$,
which is impossible. Therefore the assumption above that
all $|\lambda_i|$ are in $\{2,4,5,6\}$ was wrong.

In the cases $G_2,F_4,E_6,E_7$ the type of the subroot lattices $L_j$ 
with $[L:L_j]=3$ is unique, in the case $E_8$ there are
two possibilities, 
\begin{eqnarray*}
\begin{array}{l|l|l|l|l|l}
L & G_2 & F_4 & E_6 & E_7 & E_8\\ \hline
L_i & A_2 & A_2+\www A_2 & 3A_2 & A_5+A_2 &  E_6+A_2 ,\ A_8
\end{array}
\end{eqnarray*}

By renumbering the roots, we can assume $[L:L_{n+1}]=3$.

\medskip
The case ${\bf G_2}$:
The roots $\alpha_1$ and $\alpha_2$ generate an $A_2$
lattice and thus are long. Therefore $\alpha_3$ is short.
At least one of $\alpha_1$ and $\alpha_2$ is not orthogonal 
to $\alpha_3$. That root and $\alpha_3$ form a $\Z$-basis of $L$.

\medskip
The cases ${\bf F_4,E_6,E_7}$ and the case 
${\bf (E_8,E_6+A_2)}$:
The sublattice $L_{n+1}\subset L$ (with $\Z$-basis 
$\alpha_1,...,\alpha_n$) contains one orthogonal summand
$\www L_1$ of type $A_2$. Suppose that $\alpha_1$ and $\alpha_2$
form a $\Z$-basis of this lattice $\www L_1$.
Then $\www L_2:=\sum_{i\in\{1,2,n+1\}}\Z\alpha_i$ is a subroot
lattice with $\www L_1\subset \www L_2\subset L$ and 
$\www L_2\not\subset L_{n+1}$
and $L_{n+1}\not\subset \www L_2$.
Because of $0=\sum_{i=1}^{n+1}\lambda_i\alpha_i$ and
all $\lambda_i\neq 0$ and $n+1>3$, the root $\alpha_{n+1}$
is not in $\www L_{1,\Q}$, so the lattice $\www L_2$ has rank 3.
The sum 
$$\sum_{i\in\{1,2,n+1\}}\lambda_i\alpha_i
=-\sum_{i=3}^n\lambda_i\alpha_i$$ 
is in the sum
$\www L_{1,\R}^\perp\cap L_{n+1}$
of the other orthogonal summands of $L_{n+1}$
and in the rank one $\Z$-lattice
$\www L_{1,\R}^\perp\cap \www L_2$.
In fact, it is a generator of this rank one lattice:
This is equivalent to $c_1=1$ where
$$c_1:=\gcd(\lambda_1,\lambda_2,\lambda_{n+1}).$$
If $c_1>1$, then $c_1^{-1}\sum_{i\in\{1,2,n+1\}}\lambda_i\alpha_i$
were in the root lattice $\www L_2$. But then also
$c^{-1}_1\sum_{i=3}^n\lambda_i\alpha_i$ were in the root lattice
$L_{n+1}$, thus $c$ would divide $\lambda_3,...,\lambda_n$.
Then $\gcd(\lambda_1,...,\lambda_{n+1})>1$ which is not true.

The root lattice $\www L_1$ is either of type $A_3$ or
of type $B_3$. In both cases, the following lemma gives
the claim.

\begin{lemma}\label{t4.3}
In both cases, at least one of $\lambda_1$ and $\lambda_2$
is equal to $\pm 1$.
\end{lemma}

{\bf Proof:}
(a) The case $\www L_2$ of type $A_3$:
Embed $\www L_{2,\R}$ as usual into a Euclidean space
$\R_4$ with ON-basis $\www e_1,\www e_2,\www e_3,\www e_4$ such that 
$\Phi(\www L_2)=\{\pm (\www e_i-\www e_j)\, |\, 1\leq i<j\leq 4\}$
and $\Phi(\www L_1)=\{\pm (\www e_i-\www e_j), |\, 1\leq i<j\leq 3\}$.
A generator of the rank one $\Z$-lattice 
$\www L_{1,\R}^\perp \cap\www L_2$
is obviously $\www e_1+\www e_2+\www e_3-3\www e_4$. Thus
$$\www e_1+\www e_2+\www e_3-3\www e_4
=\pm\sum_{i\in\{1,2,n+1\}}\lambda_i\alpha_i.$$
Lemma \ref{t4.2} (a) applies to the graph
$\GG(\{\alpha_1,\alpha_2,\alpha_{n+1}\})$.
The graph is a tree, as
$\alpha_1,\alpha_2,\alpha_{n+1}$ is a $\Z$-basis of $\www L_2$.
At least two of the four vertices $\www e_1,\www e_2,\www e_3,\www e_4$
are leaves, 
so at least one of the three vertices $\www e_1,\www e_2,\www e_3$ is a leaf.
Let $\alpha_j$ with $j\in\{1,2\}$ give the edge which contains
this vertex. The coefficient $1$ of this vertex must be equal
to $\pm\lambda_j$ because the other two terms in the sum
$\pm\sum_{i\in\{1,2,n+1\}}\lambda_i\alpha_i$ have no contribution
to the coefficient of this vertex.

\medskip
(b) The case $\www L_2$ of type $B_3$:
Embed $\www L_{2,\R}$ as usual into a Euclidean space
$\R_3$ with ON-basis $\www e_1,\www e_2,\www e_3$ such that 
$\Phi(\www L_2)=\{\pm \www e_i\pm \www e_j)\}\cup\{\pm e_i\}$
and $\Phi(\www L_1)=\{\pm (\www e_i-\www e_j), |\, 1\leq i<j\leq 3\}$.
A generator of the rank one $\Z$-lattice 
$\www L_{1,\R}^\perp\cap \www L_2$
is obviously $\www e_1+\www e_2+\www e_3$. Thus
\begin{eqnarray}\label{4.3}
\www e_1+\www e_2+\www e_3=\pm\sum_{i\in\{1,2,n+1\}}\lambda_i\alpha_i.
\end{eqnarray}

Lemma \ref{t4.2} applies to the graph
$\GG(\{\alpha_1,\alpha_2,\alpha_{n+1}\})$.
As $\alpha_1$ and $\alpha_2$ are long roots, the graph
is a tree with one marked vertex.
Then at least one vertex $\www e_j$ of the 3 vertices $\www e_1,\www e_2,\www e_3$ 
has no marking and is a leaf. Then $\lambda_j=\pm 1$ for 
the same reason as in (a).
\hfill$\Box$

\medskip
The case ${\bf (E_8,A_8)}$:
We can choose a root basis $\delta_1,...,\delta_8$ of $L$
and an additional root $\delta_9$ such that they give
rise to the extended Dynkin diagram in remark \ref{t2.4} (v)
and such that the subroot lattice $L_{n+1}$ is generated by
$\delta_1,\delta_2,\delta_4,...,\delta_9$.
Further, we can embed $L_\R=L_{n+1,\R}$ into a Euclidean
space $\R^9$ with ON-basis $\www e_1,...,\www e_9$ such that
\begin{eqnarray}\label{4.4}
(\delta_1,\delta_2,\delta_4,...,\delta_9)=
(\www e_1-\www e_2,\www e_2-\www e_3,\www e_3-\www e_4,...,
\www e_8-\www e_9).
\end{eqnarray}
Part (a) of the following lemma tells how $\Phi(L)$
can be expressed using the $\www e_1,...,\www e_9$.
Part (b) solves the case $(E_8,A_8)$ and finishes the proof of
theorem \ref{t4.1}.

\begin{lemma}\label{t4.4}
(a) 
\begin{eqnarray}\label{4.5}
\delta_3 &=& \frac{1}{3}
\left(-2\sum_{i=1}^3\www e_i + \sum_{i=4}^9\www e_i\right),\\
\Phi(A_8)&=&\{\pm (\www e_i-\www e_j)\, 
|\, 1\leq i<j\leq 8\},\nonumber\\
\Phi(E_8)&=& \Phi(A_8)\cup
\{\pm\frac{1}{3}\left(-2\sum_{i\in I_1}\www e_i + 
\sum_{i\in I_2}\www e_i\right) \, 
|\, I_1\cup I_2=\{1,...,9\},\nonumber\\ 
&&\hspace*{5.5cm} |I_1|=3,|I_2|=6\}.\label{4.6}
\end{eqnarray}
(b) Above, at least one of the $\lambda_j$
with $1\leq j\leq 8$ is equal to $\pm 1$.
\end{lemma}

{\bf Proof:}
(a) \eqref{4.5} follows from \eqref{4.4}
and the relation 
$$0=2\delta_1+4\delta_2+3\delta_3+6\delta_4+5\delta_5
+4\delta_6+3\delta_7+2\delta_8+\delta_9,$$
see lemma \ref{t3.7} (c).
As $\Phi(E_8)$ contains $\delta_3$, it contains the 
combination of $\www e_i$ on the right hand side of \eqref{4.5}.
As the Weyl group $W(A_8)\subset W(E_8)$ consists of all
permutations of $\www e_1,...,\www e_9$, the root system
$\Phi(E_8)$ contains the right hand side of \eqref{4.6}.
Counting the size of the right hand side, one finds
$$2\begin{pmatrix}9\\2\end{pmatrix}
+2\begin{pmatrix}9\\3\end{pmatrix}=72+168=240=|\Phi(E_8)|,$$
thus equality holds in \eqref{4.6}.

\medskip
(b) The roots $\alpha_1,...,\alpha_8$ form a $\Z$-basis of
$L_{n+1}$. The root $\alpha_9$ must be a root in 
$\Phi(L)-\Phi(L_{n+1})$, so it must be
$$\alpha_9=\pm\frac{1}{3}
\left(-2\sum_{i=1}^3\www e_{\pi(i)} + 
\sum_{i=4}^9\www e_{\pi(i)}\right)\quad
\textup{for some }\pi\in S_9.$$
The graph $\GG(\{\alpha_1,...,\alpha_8\})$ is a tree 
by lemma \ref{t4.2} (a).

\medskip
{\bf 1st case,} at least one of the roots $\www e_{\pi(4)},...,
\www e_{\pi(9)}$ is a leaf in this graph:
Let $\alpha_j$ be the only edge which contains this leaf.
Then $\lambda_j\alpha_j$ contains the only contribution to
the leaf, in the right hand side of the following formula,
\begin{eqnarray*}
-2\sum_{i=1}^3\www e_{\pi(i)} + 
\sum_{i=4}^9\www e_{\pi(i)}
=\pm 3\alpha_9 = \pm\sum_{j=1}^8 \lambda_j\alpha_j.
\end{eqnarray*}
Thus then $\lambda_j=\pm 1.$

\medskip
{\bf 2nd case,} none of the roots $\www e_{\pi(4)},...,
\www e_{\pi(9)}$ is a leaf in the graph:
Then there are two or three leafs, and they form a subset
of the set $\{\www e_{\pi(1)},\www e_{\pi(2)},\www e_{\pi(3)}\}$.
For one of these leafs, the number of vertices
on the path from this leaf to the branching vertex (in the 
case of three leafs) or to the unique inner vertex 
which is in $\{\www e_{\pi(1)},\www e_{\pi(2)},\www e_{\pi(3)}\}$,
is maximal. Then the first two edges within this path,
which starts at the leaf, have the coefficients $\lambda_j$
with values $\pm 2,\pm 1$. So $\pm 1$ arises. \hfill$\Box$

\medskip
This finishes the proof of theorem \ref{t4.1}\hfill$\Box$

\section{Reduced presentations of Weyl group elements}\label{s5}
\setcounter{equation}{0}

Carter studied and classified the conjugacy classes of the elements
of the Weyl groups of the irreducible root lattices.
Here we will review a part of his results and extend them.
Crucial are the (in the inhomogeneous cases new) notions
of quasi Coxeter elements and strict quasi Coxeter elements.
The control of these elements reduces the classification of conjugacy
classes of Weyl group elements to the control of subroot lattices 
in section \ref{s3}. But first some definitions will be given.

\begin{definition}\label{t5.1}
Let $(L,(.,.),\Phi)$ be a p.n. root lattice with Weyl group $W$.

\medskip
(a) For any element $w\in W$ any tuple $(\alpha_1,...,\alpha_k)\in \Phi^k$
with $k\in\Z_{\geq 0}$ and 
\begin{eqnarray}\label{5.1}
w=s_{\alpha_1}\circ...\circ s_{\alpha_k}
\end{eqnarray}
is a {\it presentation of $w$}. Its {\it length} is $k\in\Z_{\geq 0}$.
The length $l(w)\in\Z_{\geq 0}$ of $w$ is the minimum of the
lengths of all presentations. A presentation with $k=l(w)$ is called
{\it reduced}. 
The subroot lattice of a presentation $(\alpha_1,...,\alpha_k)$ 
is $L_1:=\sum_{i=1}^k\Z\cdot\alpha$. 
The {\it index} of the presentation is the index
$[L\cap L_{1,\Q}:L_1]\in\Z_{\geq 1}$ of the subroot lattice $L_1$.

\medskip
(b) An element $w$ is of {\it maximal length} if $l(w)=n:=$
the rank of the root lattice. 

\medskip
(c) For any element $w\in W$ and any $\lambda\in S^1$ define
\begin{eqnarray}\label{5.2}
V_\lambda(w)&:=& \ker(w-\lambda\cdot\id)\subset L_\C,\\
V_{\neq 1}(w)&:=&\bigoplus_{\lambda\neq 1}V_\lambda(w)
\supset V_{\neq 1,\R}:=L_\R\cap V_{\neq 1}(w),\label{5.3}
\end{eqnarray}
and analogously $V_{\neq 1,\Q}(w), V_{\neq 1,\Z}(w)$
$V_{1,\R}(w), V_{1,\Q}(w),V_{1,\Z}(w)$.
Of course $V_{\neq 1,\R}=V_{1,\R}^\perp$.
\end{definition}

\begin{lemma}\label{t5.2}
Let $(L,(.,.),\Phi)$ be a p.n. root lattice with Weyl group $W$.

\medskip
(a) \cite[Lemmata 2 and 3]{Ca}
A presentation $(\alpha_1,...,\alpha_k)$ 
of an element $w\in W$ is reduced if and only if
$\alpha_1,...,\alpha_k$ are linearly independent (in $L_\Q$).
The subroot lattice $L_1\subset L$ of a reduced presentation satisfies
\begin{eqnarray}\label{5.4}
\left(\bigoplus_{i=1}^{l(w)}\Q\cdot\alpha_i =\right)L_{1,\Q} &=& 
V_{\neq 1,\Q}(w)\\
\textup{and especially}\qquad l(w) &=& \dim V_{\neq 1,\Q}(w).\label{5.5}
\end{eqnarray}
So, the subroot lattices of all reduced presentations of $w$ generate
the same subspace of $L_\C$, and it is $V_{\neq 1}(w)$.

\medskip
(b) \cite[Satz 3.2]{Kl}\cite[Satz 3.2.3]{Vo}
If $(L,(.,.),\Phi)$ is a homogeneous root lattice,
then all reduced presentations of one element $w\in W$ have the same index.
\end{lemma}

The following definition of a quasi Coxeter element 
is in the homogeneous cases due to Voigt
\cite[Def. 3.2.1]{Vo} and in the inhomogeneous cases new.

\begin{definition}\label{t5.3}
Let $(L,(.,.),\Phi)$ be a p.n. root lattice of rank $n\in\Z_{>0}$ 
with Weyl group $W$.

\medskip
(a) An element $w\in W$ 
is a {\it quasi Coxeter element}
if a reduced presentation of $w$ exists whose
subroot lattice is the full root lattice $L$.
Of course then it is of {\it maximal length} $l(w)=n$.

\medskip
(b) An element $w\in W$
is a {\it strict quasi Coxeter element}
if the subroot lattice of any reduced presentation 
is the full root lattice $L$.
Of course then it is a quasi Coxeter element.
\end{definition}

\begin{remarks}\label{t5.4}
(i) An element $w$ in the Weyl group of a p.n. root lattice
has many presentations. Often there are several presentations such
that the isomorphisms classes of their subroot lattices are different.
In the homogeneous cases at least their indices are equal.
But in the inhomogeneous cases, even their indices can differ.

\medskip
(ii) In a homogeneous root lattice, lemma \ref{t5.2} (b) implies
that there the notions quasi Coxeter element and strict quasi Coxeter element
coincide. But in any irreducible inhomogeneous root lattice, there are quasi
Coxeter elements which are not strict quasi Coxeter elements. See theorem
\ref{t5.6}. 

\medskip
(iii) Of course, if $(\alpha_1,...,\alpha_k)$ is a reduced presentation of a
Weyl group element $w$, then $w$ is a quasi Coxeter element in the
subroot lattice $L_1$ of this presentation. And of course, any 
Weyl group element has a reduced presentation such that it is a strict quasi Coxeter
element in the subroot lattice $L_1$ of this presentation.

\medskip
(iv) Let $L=\bigoplus_{k\in K}L_k$ be the decomposition of
a p.n. root lattice into an orthogonal sum of irreducible p.n.  root lattices,
and let $w\in W$ be a [strict] quasi Coxeter element.
Then it decomposes into a product $\prod_{k\in K}w_k$ of commuting
elements $w\in W(L_k)$, and $w_k$ is a [strict] quasi Coxeter element
in $L_k$. 

\medskip
(v) Recall that a Coxeter element in an irreducible root lattice
is an element $w\in W$ which has a presentation $(\alpha_1,...,\alpha_n)$
such that $\alpha_1,...,\alpha_n$ form a root basis.
Because their Dynkin diagram is a tree, lemma 1 in \cite[Ch. V \S 6]{Bo}
implies that the products of $s_{\alpha_1},...,s_{\alpha_n}$ in any
order are conjugate. As all root bases are conjugate, all
Coxeter elements are conjugate.
Obviously the Coxeter elements are quasi Coxeter elements.
It turns out that they are even strict quasi Coxeter elements, see
theorem \ref{t5.6}.

\medskip
(vi) Carter's work \cite{Ca} on the classification of Weyl group
elements gives in a direct way the classification of the quasi Coxeter
elements in the irreducible homogeneous root lattices
and in a less direct way the classification of the strict 
quasi Coxeter elements in the irreducible inhomogeneous root lattices.
In theorem \ref{t5.6}  
these classifications will be given, and also the classification
of the quasi Coxeter elements in the irreducible inhomogeneous root lattices. 

\medskip
(vii) Recall the description of the Weyl group $W$ in remark \ref{t2.4} (iii)
for the root lattices of the types $A_n,B_n,C_n,D_n$ in theorem \ref{t2.3}
(d): $W(A_n)\cong S_{n+1}$, $ W(B_n)=W(C_n)\cong \{\pm 1\}^n\rtimes S_n$.
A signed permutation in $\{\pm 1\}^n\rtimes S_n$
will be called {\it positive} if the number of
sign changes in it is even, it will be called {\it negative} if the number of
sign changes in it is odd. The subgroup $W(D_n)\subset W(B_n)=W(C_n)$ 
consists of the positive signed permutations.

A signed cycle will be written as
$(\varepsilon_1 a_1\, \varepsilon_2 a_2 ... \varepsilon_k a_k)$ with $k\geq 1$ 
and $\varepsilon_1,...,\varepsilon_k\in\{\pm 1\}$, $a_1,...,a_k\in\{1,...,n\}$
with $a_i\neq a_j$ for $i\neq j$. It maps $\pm a_i$ to 
$\pm \varepsilon_{i+1}a_{i+1}$ for $1\leq i\leq k-1$ and $\pm a_k$ to 
$\pm \varepsilon_1 a_1$. It is positive if $\prod_j\varepsilon_j=1$
and negative if $\prod_j\varepsilon_j=-1$. 
Its {\it support} is defined to be $\{a_1,...,a_k\}$. 

Any signed permutation is up to the order
a unique product of signed cycles (=cyclic permutations) such that their
supports are disjoint and the union of the supports is $\{1,...,n\}$. 
They are called the signed cycles of the permutation. 
Here cycles of length one are used.
For example $\id=(1)(2)...(n)$ and $-\id =(-1)(-2)...(-n)$.
\end{remarks}

\begin{remarks}\label{t5.5}
(i) 
Carter classified in \cite{Ca} the conjugacy classes of Weyl group
elements for all irreducible root lattices. A crucial point was the
proof that any element $w$ can be written as a product $w=w_1w_2$
where $w_1$ and $w_2$ are involutions with $V_{-1}(w_1)\cap V_{-1}(w_2)=\{0\}$
(proposition 38 and corollary (ii) in \cite{Ca}). 

By \cite[lemma 5]{Ca}, any involution has a reduced 
presentation $(\alpha_1,...,\alpha_k)$ which 
consists of pairwise orthogonal roots.
The composition of two such reduced presentations of two involutions 
$w_1$ and $w_2$ with $V_{-1}(w_1)\cap V_{-1}(w_2)=\{0\}$ 
is a reduced presentation of $w=w_1w_2$. 
Its generalized Dynkin diagram is a graph whose cycles 
(if any exist) have all even length. 
In \cite[theorem A]{Ca} all graphs are classified which
have the following properties: 
The graph contains cycles, all cycles have even length,
the graph is a generalized Dynkin diagram of a presentation of an element $w=w_1w_2$
with $w_1$ and $w_2$ as above, the subroot lattice of the presentation
is the full lattice, and $w$ is not contained in the Weyl group of a 
subroot system. The graphs are labelled
$D_n(a_k),E_6(a_k),E_7(a_k),E_8(a_k),F_4(a_1)$ and also $D_n(b_{n/2-1})$
if $n$ is even and $E_7(b_2),E_8(b_3),E_8(b_5)$. 

In fact, the graphs in \cite{Ca} are simplified by not distinguishing
normal and dotted edges. The generalized Dynkin diagrams are obtained from
the graphs in \cite{Ca} by replacing some edges by dotted edges such that
any cycle obtains an odd number of dotted edges. This is possible.

It turns out that the graphs $D_n(a_k),E_6(a_k),E_7(a_k),E_8(a_k),F_4(a_1)$
correspond to conjugacy classes of Weyl group elements,
and that these include the elements with presentations giving rise to the
graphs $D_n(b_{n/2-1})$ ($n$ even) and $E_7(b_2),E_8(b_3),E_8(b_5)$.
These Weyl group elements are strict quasi Coxeter elements,
because they are not contained in the Weyl group of a subroot lattice.
They are not Coxeter elements \cite{Ca}.
They and the Coxeter elements 
are the only strict quasi Coxeter elements (theorem \ref{t5.6} below).

\medskip
(ii)
Recall that the Coxeter elements in $W(A_n)$ are the cycles of length
$n+1$ in $S_{n+1}$, 
the Coxeter elements in $W(B_n)=W(C_n)=W(BC_n)$
are the negative cycles of length $n$ in $\{\pm 1\}\rtimes S_n$,
and the Coxeter elements in $W(D_n)$ are the products of two
negative cycles of lengths 1 and $n-1$. The products of two negative
cycles of lengths $k$ and $n-k$ for $2\leq k\leq [n/2]$ form
the conjugacy class $D_n(a_{k-1})$ in $W(D_n)$.

In the second column in the tables 5.1 and 5.2,
$A_n,B_n,C_n,D_n,E_6,E_7,E_8,F_4$ and $G_2$
denote the conjugacy classes of the Coxeter elements.
The root lattice of type $F_4$ contains subroot lattices of types
$B_4$ and $D_4$. In $W(F_4)$ the symbols 
$B_4,C_3+A_1,D_4(a_1)$ denote the conjugacy classes in $W(F_4)$ of the
Coxeter elements in $W(B_4)$ and $C_3+A_1$ and of the quasi Coxeter
elements of type $D_4(a_1)$ in $W(D_4)$.
The Coxeter elements of the subroot lattice of type $A_2$ in $G_2$
give rise to a conjugacy class in $W(G_2)$ denoted by $A_2$.
\end{remarks}

Theorem \ref{t5.6} gives the classification of the quasi Coxeter
elements and the strict quasi Coxeter elements for the irreducible
p.n. root lattices. A good part of it is due to \cite{Ca}.

\begin{theorem}\label{t5.6}
Let $(L,(.,,.),\Phi)$ be one of the irreducible p.n. root lattices 
in theorem \ref{t2.3} (d). 
The following tables 5.1 and 5.2 list the conjugacy classes of the 
strict quasi Coxeter elements and
in the inhomogeneous cases the conjugacy classes of the
quasi Coxeter elements. See the remarks \ref{t5.5} for
the notations.

\medskip

{\bf Table 5.1:}
\begin{eqnarray*}
\begin{array}{l|l}
 & \textup{strict quasi Coxeter el. = quasi Coxeter el.}\\ \hline
A_n & A_n \\
D_n & D_n, D_n(a_1),...,D_n(a_{[n/2-1]})\\
E_6 & E_6, E_6(a_1),E_6(a_2)\\
E_7 & E_7, E_7(a_1),E_7(a_2),E_7(a_3),E_7(a_4)\\
E_8 & E_8, E_8(a_1),...,E_8(a_8)
\end{array}
\end{eqnarray*}

{\bf Table 5.2:}
\begin{eqnarray*}
\begin{array}{l|l|l}
 & \textup{strict quasi Coxeter el.} & \textup{quasi Coxeter el.}\\ \hline
B_n & B_n & \textup{products of negative cycles}\\
BC_n & - & \textup{products of negative cycles}\\
C_n & C_n & C_n,D_n, D_n(a_1),...,D_n(a_{[n/2-1]})\\
F_4 & F_4,F_4(a_1) & F_4,F_4(a_1), B_4,C_3+A_1,D_4(a_1)\\
G_2 & G_2 & G_2,A_2
\end{array}
\end{eqnarray*}
Thus the quasi Coxeter elements of $C_n$ consist of the
products of one or two negative cycles.
\end{theorem}

{\bf Proof:}
It is well known that the Coxeter elements are not elements of some
proper Weyl subgroup. Therefore they are strict quasi Coxeter elements.
The other elements listed in the second columns are strict quasi Coxeter
elements because of the results of Carter \cite{Ca} discussed
in the remarks \ref{t5.5} (i).

By the same results, any other element $w\in W$
is in some proper Weyl subgroup. In the homogeneous cases, a proper
Weyl subgroup is the Weyl group of a proper subroot lattice. Therefore
then $w$ is not a strict quasi Coxeter element. This completes the proof
of table 5.1.

In the inhomogeneous cases, the fact that the second column of table 5.2
lists all strict quasi Coxeter elements 
is a consequence of the third column of table 5.2, 
in the following way. In the cases of the
root lattices of types $F_4$ and $G_2$ it is obvious that the
quasi Coxeter elements of types $B_4,C_3+A_1,D_4(a_1)$ 
and $A_2$ are not strict quasi Coxeter elements. 
In the cases of the p.n. root lattices $B_n,C_n$ and $BC_n$, 
observe
\begin{eqnarray}\label{5.6}
s_{e_i}s_{e_j} &=& s_{e_i+e_j}s_{e_i-e_j}
=s_{2e_i}s_{2e_j}\sim (-i)(-j) \qquad\textup{and}\\ \nonumber
\Z e_i+\Z e_j &\supsetneqq& \Z(e_i+e_j)+\Z(e_i-e_j)
\supsetneqq \Z 2e_i+\Z 2e_j\qquad\textup{for }i\neq j.
\end{eqnarray}
This shows that all permutations whose signed cycles contain at least
two negative cycles are not strict quasi Coxeter elements.
$BC_n$ has no strict quasi Coxeter elements because of $s_{e_i}=s_{2e_i}$.
Therefore the only elements in the third column of table 5.2 which are 
strict quasi Coxeter elements are those in the second column.

It rests to prove the third column of table 5.2.

The root lattice of type ${\bf C_n}$:
Because of $L(C_n)=L(D_n)$ the quasi Coxeter elements of
$D_n$ are also quasi Coxeter elements of $C_n$.
Let $w\in W(C_n)$ be a quasi Coxeter element of $C_n$
which is not a quasi Coxeter element of $D_n$. Then
it has a presentation $(\alpha_1,...,\alpha_n)$ such
that its subroot lattice is the full root lattice $L(C_n)$
and $A:=\{\alpha_1,...,\alpha_n\}\not\subset \Phi(D_n)$.
By lemma \ref{t4.2} (c)(ii) and (e)(ii) then the graph $\GG(A)$
is a tree and contains exactly one marked vertex.
Then $w$ is a negative cycle, so a Coxeter element of $C_n$.

The p.n. root lattices of types ${\bf B_n}$ and ${\bf BC_n}$:
Let $w\in W(B_n)=W(BC_n)$ be a quasi Coxeter element of $B_n$
or $BC_n$, and let $(\alpha_1,...,\alpha_n)$ be a 
presentation whose subroot lattice is the full root lattice
$L(B_n)=L(BC_n)$. By lemma \ref{t4.2} (b)(ii) and (d)(ii)
then the graph $\GG(A)$ for $A=\{\alpha_1,...,\alpha_n)$
is a union of trees which have each exactly one marking
and which is short. Thus $w$ is a product of negative cycles.
Vice versa, any product $w$ of negative cycles has
a presentation $(\alpha_1,...,\alpha_n)$ such that the
graph $\GG(\{\alpha_1,...,\alpha_n\})$ is a union of trees
which have each exactly one marking and which is short.
Thus $w$ is a quasi Coxeter element.

The root lattice of type ${\bf G_2}$:
Obviously its quasi Coxeter elements are the products
$s_\alpha s_\beta$ with $\alpha$ short and $\beta$ long
and $\alpha\not\perp\beta$ and the products 
$s_{\alpha_1}s_{\alpha_2}$ with $\alpha_1$ and $\alpha_2$ short
and $\alpha_2\neq\pm\alpha_1$. The elements of the first
type are the Coxeter elements of $G_2$, the elements of the
second type can also be written as products
$s_{\beta_1}s_{\beta_2}$ with $\beta_1$ and $\beta_2$
long roots and $\beta_2\neq\pm\beta_1$. They 
are the Coxeter elements of the subroot system of long roots,
which is of type $A_2$.

The root lattice of type ${\bf F_4}$:
See lemma \ref{t5.7} (b). The restriction there that
in the presentation $(\alpha_1,\alpha_2,\alpha_3,\alpha_4)$
first the short roots come and then the long roots,
is not serious. One can obtain a presentation with this
property from an arbitrary presentation using \eqref{2.2b}.
\hfill$\Box$

\begin{lemma}\label{t5.7}
Let $(L,(.,.),\Phi)$ be the root lattice of type $F_4$
in theorem \ref{t2.3} (d). 
Obviously the short roots form a root system of type $D_4$,
which is called $\www D_4$,
and the long roots form a root system of type $D_4$,
which is called $D_4$.

\medskip
(a) Let $A=\{\alpha_1,\alpha_2,\alpha_3,\alpha_4\}\subset 
\Phi(F_4)$ be a $\Z$-basis of $L(F_4)$ such that first
the short roots come and then the long roots. Then
one of the following cases holds.
\begin{list}{}{}
\item[(i)]
All four roots are short. 
\item[(ii)]
$\alpha_1,\alpha_2$ and $\alpha_3$ are short and $\alpha_4$
is long.
Then an element $w\in W(F_4)$ exist such that
$w(\alpha_1),w(\alpha_2),w(\alpha_3)$ generate the subroot
system of type $\www A_3$ which is also generated by 
$e_1,\frac{1}{2}\sum_{i=1}^4 e_i,e_2$, and then
$w(\alpha_4)=\pm e_i\pm e_j$ with $i\in\{1,2\}$
and $j\in\{3,4\}$.
\item[(iii)]
$\alpha_1$ and $\alpha_2$ are short and $\alpha_3$ and
$\alpha_4$ are long. Then $\langle\alpha_1,\alpha_2\rangle=\pm 1$
and $\langle\alpha_3,\alpha_4\rangle=\pm 1$ and 
$\R\alpha_1+\R\alpha_2\not\perp \R\alpha_3+\R\alpha_4$.
\end{list}

\medskip
(b) Let $w\in W(F_4)$ be a quasi Coxeter element,
and let $(\alpha_1,...,\alpha_4)$ be a presentation of $w$ whose
subroot lattice is the full lattice $L(F_4)$ and 
such that first the short roots come and then the long roots.
Then the cases in (a) hold, and $w$ is in each case as follows.
\begin{list}{}{}
\item[(i)]
$w$ is in $W(C_3+A_1)$ and is a Coxeter element there,
or it is in $W(D_4)$ and is a quasi Coxeter element of type
$D_4(a_1)$ there.
\item[(ii)]
$w$ is in $W(L_1)$ for some subroot lattice $L_1$ of type $B_4$
and is a Coxeter element there.
\item[(iii)]
$w$ is a Coxeter element of type $F_4$ or a strict quasi
Coxeter element of type $F_4(a_1)$.
\end{list}
\end{lemma}

{\bf Proof:}
(a) The following obvious statements will be used:

\begin{list}{}{}
\item[(A)]
The root lattice $L(F_4)$ contains the subroot lattice 
$L_2$ of type $B_4$ with root system
$$\Phi(L_2)=\{\pm e_i\, |\, i\in\{1,2,3,4\}\}
\cup \{\pm e_i\pm e_j\, |\, 1\leq i<j\leq 4\}.$$
It has the same long roots as $\Phi(F_4)$, but less short roots.
\item[(B)]
For any short root $\beta_1$
a Weyl group element $w$ exists such that $w(\beta_1)=e_1$.
If $\beta_2$ is a short root with $\beta_1\perp\beta_2$
then $w(\beta_2)\in\{\pm e_2,\pm e_3,\pm e_4\}$, and
$w$ can be chosen such that $w(\beta_2)=e_2$.
\end{list}

The case that all four roots $\alpha_1,...,\alpha_4$ are long
is impossible because they would only generate the subroot
lattice $L(D_4)$. The case that the three roots 
$\alpha_1,\alpha_2,\alpha_3$ are long and $\alpha_4$ is short,
is also impossible, because by (B) a Weyl group element
exists such that $w(\alpha_4)=e_1$, and then all four images 
$w(\alpha_i)$ are in $L_2$.
Thus either two roots are short and two roots are long,
or three roots are short and one root is long,
or all four roots are short.

Consider the case that $\alpha_1,\alpha_2,\alpha_3$ are short
and $\alpha_4$ is long. Then $\alpha_1,\alpha_2$ and $\alpha_3$
generate a subroot system of rank 3 of $\Phi(\www D_4)$.
Only the two types $\www A_3$ and $3\www A_1$ are possible a priori.
Here the type $3\www A_1$ is not possible, because then by (B) an
element $w\in W(F_4)$ exists such that $w(\alpha_i)=e_i$
for $i\in\{1,2,3\}$, and these roots and any
long root $w(\alpha_4)$ are in $L_2$.
By theorem \ref{t3.3} (b), any two subroot systems of type
$\www A_3$ of the subroot system $\www D_4$ are conjugate
by an element of $W(\www D_4)$. This shows the first
half of part (ii). Obviously $w(\alpha_4)=\pm e_i\pm e_j$
with $i\in\{1,2\}$ and $j\in\{3,4\}$. This gives part (ii).

Consider the case that $\alpha_1$ and $\alpha_2$ are short and
$\alpha_3$ and $\alpha_4$ are long.
If $\alpha_1\perp\alpha_2$ then by (B) 
$\alpha_1,...,\alpha_4$ are mapped by a suitable element
$w\in W(F_4)$ into the subroot lattice $L_2$.
Therefore $\langle\alpha_1,\alpha_2\rangle=\pm 1$.
Furthermore 
$$\R\alpha_1+\R\alpha_2\not\perp\R\alpha_3+\R\alpha_4$$
because else the four roots would generate a reducible
subroot lattice. An element $w\in W(F_4)$ exists such
that 
$$w(\alpha_1)=\pm e_1,w(\alpha_2)=\frac{1}{2}\sum_{i=1}^4 e_i.$$
If $\alpha_3\perp\alpha_4$ then either 
$$w(\alpha_3)=\varepsilon_1(e_i+\varepsilon_2 e_j)\textup{ and }
w(\alpha_4)=\varepsilon_3(e_i-\varepsilon_2 e_j)$$
for some $\varepsilon_1,\varepsilon_2,\varepsilon_3\in\{\pm 1\}$
and some $i,j$ with $1\leq i<j\leq 4$, or 
$$w(\alpha_3)=\varepsilon_1(e_i+\varepsilon_2 e_j)\textup{ and }
w(\alpha_4)=\varepsilon_3(e_k-\varepsilon_4 e_l)$$
for some $\varepsilon_1,\varepsilon_2,\varepsilon_3,
\varepsilon_4\in\{\pm 1\}$
and some $i,j,k,l$ with $\{i,j,k,l\}=\{1,2,3,4\}$.
One sees easily with some case discussion 
that in both cases $w(\alpha_1),...,w(\alpha_4)$
do not generate $L(F_4)$.

\medskip
(b) Of course, the cases in (a) hold.

{\bf The case (i):} $w$ is in $W(\www D_4)$ and is either a 
Coxeter element there or a quasi Coxeter element of type
$\www D_4(a_1)$. In the first case, $w$ is conjugate to
$$s_{\frac{1}{2}\sum_ie_i}s_{e_1}s_{e_2}s_{e_3}
=s_{\frac{1}{2}\sum_ie_i}s_{e_1}s_{e_2+e_3}s_{e_2-e_3},$$
which is in $W(C_3+A_1)$ and which is a Coxeter element there.
In the second case, $w$ is conjugate to
$$s_{\frac{1}{2}\sum_ie_i}s_{\frac{1}{2}(e_1+e_2-e_3-e_4)}
s_{e_1}s_{e_3}
=s_{e_1+e_2}s_{e_3+e_4}s_{e_1+e_3}s_{e_1-e_3},$$
which is in $W(D_4)$ and which is a quasi Coxeter element
of type $D_4(a_1)$ there.

{\bf The case (ii):} The element $w$ is conjugate to
$$s_{\frac{1}{2}\sum_ae_a}s_{e_1}s_{e_2}s_{\beta}
=s_{\frac{1}{2}\sum_ae_a}s_{e_1-e_2}s_{e_1+e_2}s_{\beta}$$
for some $\beta=e_i+\varepsilon e_k$ with $i\in\{1,2\}$
and $k\in\{3,4\}$ and $\varepsilon\in\{\pm 1\}$.
This is in $W(L_3)$ for a subroot lattice $L_3$ of type $B_4$.
In the case $\varepsilon=-1$ the generalized Dynkin 
diagram of the four roots on
the right hand side is (up to the distinction between dotted 
and normal edges) the $B_4$ Dynkin diagram, so then
the element is a Coxeter element in $W(L_3)$.
In the case $\varepsilon=1$, the right hand side is equal to
$$s_{\frac{1}{2}\sum_ae_a}s_{e_1-e_2}s_{-e_j+e_k}s_{e_1+e_2},$$
where $j$ is determined by $\{i,j\}=\{1,2\}$. 
This is again a Coxeter element in $W(L_3)$.

{\bf The case (iii):}
Using \eqref{2.2b} for $\alpha_1$ and $\alpha_2$,
one can suppose $\R\alpha_2\not\perp\R\alpha_3+\R\alpha_4$.
After conjugation, one can suppose
$$\alpha_1=\frac{1}{2}\sum_{i=1}^4 \pm e_i,\ 
\alpha_2=\pm e_1,\ \{\alpha_3,\alpha_4\}\subset
\{\pm(e_i-e_j)\, |\, 1\leq i< j\leq 3\}.$$
Using \eqref{2.2b} for $\alpha_3,\alpha_4$ and changing possibly
some signs and conjugating possibly again, one can suppose
$$\alpha_1=\frac{1}{2}(e_1+\varepsilon_2e_2+\varepsilon_3e_3+e_4),
\ \alpha_2=e_1,\ \alpha_3=e_2-e_3,\ \alpha_4=e_1-e_2$$
for some $\varepsilon_2,\varepsilon_3\in\{\pm 1\}$.

In the case $(\varepsilon_2,\varepsilon_3)=(1,1)$, the generalized Dynkin
diagram of the roots $\alpha_1,\alpha_2,\alpha_3,\alpha_4$ is 
(up to the distinction between dotted and normal edges)
a Dynkin diagram of type $F_4$. Thus $w$ is a Coxeter element
in $W(F_4)$. 

In the case $(\varepsilon_2,\varepsilon_3)=(1,-1)$,
the element $w$ is conjugate to the product of the two involutions
$s_{\alpha_4}s_{\alpha_1}$ and $s_{\alpha_2}s_{\alpha_3}$
with admissible diagram of type $F_4(a_1)$.
Thus it is a quasi Coxeter element in $W(F_4)$ of type $F_4(a_1)$.

In the case $(\varepsilon_2,\varepsilon_3)=(-1,-1)$, the element $w$ is
$$w=s_{\alpha_1}s_{\alpha_2}s_{\alpha_3}s_{\alpha_4}
=s_{\alpha_2}s_{s_{\alpha_2}(\alpha_1)}s_{s_{\alpha_3}(\alpha_4)}s_{\alpha_3}
=s_{e_1}s_{\frac{1}{2}(-e_1-e_2-e_3+e_4)}s_{e_1-e_3}s_{e_2-e_3}.$$
The generalized Dynkin diagram of the roots on the right hand side 
is (up to the distinction between dotted and normal edges)
a Dynkin diagram of type $F_4$. Thus $w$ is a Coxeter element in $W(F_4)$.

In the case $(\varepsilon_2,\varepsilon_3)=(-1,1)$, the element $w$ is
$$w=s_{\alpha_1}s_{\alpha_2}s_{\alpha_3}s_{\alpha_4}
=s_{\alpha_1}s_{\alpha_3}s_{\alpha_2}s_{\alpha_4}
=s_{\alpha_1}s_{\alpha_3}s_{\alpha_4}s_{s_{\alpha_4}(\alpha_2)}
=s_{\alpha_1}s_{\alpha_3}s_{\alpha_4}s_{e_2}.$$
This is conjugate to the element
$$s_{e_2}s_{\alpha_1}s_{\alpha_3}s_{\alpha_4}
=s_{s_{e_2}(\alpha_1)}s_{e_2}s_{s_{\alpha_3}(\alpha_4)}s_{\alpha_3}
=s_{\frac{1}{2}\sum_i e_i}s_{e_2}s_{e_1-e_3}s_{e_2-e_3}.$$
The generalized Dynkin diagram of the roots on the right hand side 
is (up to the distinction between dotted and normal edges)
a Dynkin diagram of type $F_4$. Thus $w$ is a Coxeter element in $W(F_4)$.
\hfill$\Box$

\begin{remarks}\label{t5.8}
(i) In the tables 7--11 in \cite{Ca} all conjugacy classes
of elements of the Weyl groups of the root lattices
of types $G_2,F_4,E_6,E_7$ and $E_8$ are listed in the following
form. For any conjugacy class one element and one presentation of it 
as a strict quasi Coxeter element is chosen. The tables show the
isomorphism class of the pair of full lattice and
subroot lattice and the type of the strict quasi Coxeter element.

\medskip
(ii) Theorem \ref{t5.10} below gives more information
for the root lattices of types $G_2,F_4,E_6,E_7$ and $E_8$.
For $F_4$ it lists for any conjugacy class all 
isomorphism classes of pairs of full lattice and
subroot lattice, which turns up as subroot lattice of a
presentation as a quasi Coxeter element, and the type
of the quasi Coxeter element. This gives all types of
reduced presentations for any element. 
For $G_2,E_6,E_7$ and $E_8$ one can extract the same information 
from theorem \ref{t5.10} and the tables 3.6, 3.2, 3.3 and 3.4. 

\medskip
(iii) For the tables in theorem \ref{t5.10}, the notations
in table 5.2 have to be refined: There are three conjugacy 
classes of quasi Coxeter elements in $F_4$ 
which have also presentations as Coxeter elements 
in subroot lattices of type $B_4$ 
and $C_3+A_1$ respectively as a quasi Coxeter element 
of type $D_4(a_1)$ in the subroot lattice of type $D_4$.
The presentations of these elements as quasi Coxeter elements
in $F_4$ are now called $F_4(a_2),F_4(a_3)$ and $F_4(a_4)$.

Analogously, the presentations as quasi Coxeter elements in
$W(G_2)$ of those elements which have also presentations as
Coxeter elements in $A_2$ (the subroot lattice of long roots)
are denoted by $\www A_2$.

For $2\leq k\leq 4$, the presentations as quasi Coxeter elements
in $W(B_k)$ of those elements
in $W(B_k)$ which are products of negative cycles of
lengths $l_1,..,l_r$ with $l_1+...+l_r=k$ are denoted
by $B_k(l_1,...,l_r)$. The case $B_k(k)$ is also denoted by $B_k$.
This will be used in the table 5.4 for $F_4$.
Similarly, $C_3(2,1)$ is used there.

\medskip
(iv) The complete control in theorem \ref{t5.10}
on the reduced presentations
of the Weyl group elements of the irreducible
root lattices allows to determine the number
$k_4(L,w)$ in definition \ref{t5.9}.
\end{remarks}

\begin{definition}\label{t5.9}
Let $(L,(.,.),\Phi)$ be a p.n. root lattice, and let
$w$ be a Weyl group element. Recall $k_2(L,L_1)$ from theorem \ref{t3.8}. 
Define the number
\begin{eqnarray}\label{5.8}
k_4(L,w):=\min(k_2(L,L_1)& | & \textup{a reduced presentation of }
w\\
&& \textup{with subroot lattice }L_1\textup{ exists.}\} \nonumber
\end{eqnarray}
\end{definition}

This number will be important in section \ref{s6}.
Because of $k_1(L,L_1)=k_2(L,L_1)$ (theorem \ref{t3.8} (b)),
\begin{eqnarray}\label{5.9}
k_4(L,w)\geq \dim L_\Q/V_{\neq 1,\Q}(w)=n-l(w).
\end{eqnarray}
Equality holds if and only if a reduced presentation with
subroot lattice $L_1=V_{\neq 1,\Q}\cap L$ exists.
This is the unique primitive subroot lattice $L_1$ 
with $L_{1,\Q}=V_{\neq 1,\Q}$. Often equality holds, often not.

\begin{theorem}\label{t5.10}
Let $(L,(.,.),\Phi)$ be one of the irreducible p.n. root lattices 
in theorem \ref{t2.3} (d). 

\medskip
(a) In the cases $A_n$, $B_n$ and $BC_n$, 
\begin{eqnarray}\label{5.10}
k_4(L,w)=n-l(w).
\end{eqnarray}

\medskip
(b) Consider in the cases $C_n$ and $D_n$ a Weyl group element $w$ 
which is a product of $r$ positive cycles and $s$ negative cycles
with disjoint supports whose union is $\{1,...,n\}$
(remark \ref{t5.4} (vii)). (In the case of $D_n$, $s$ is even.) 
Then
$$r=n-l(w),$$
and then any reduced presentation with
subroot lattice $L_1$ with minimal $k_2(L,L_1)$ satisfies
\begin{eqnarray}\label{5.11}
L/L_1 &\cong& \Z^r\times \Z_2^{[(s+1)/2]},\\ 
k_4(L,w) &=& k_2(L,L_1) = n-l(w)+\left[\frac{s+1}{2}\right] .\label{5.12}
\end{eqnarray}

\medskip
(c) In the cases $G_2,E_6,E_7$ and $E_8$, for the big majority 
of the Weyl group elements there is only one type of reduced
presentations. That means, the pairs $(L,L_1)$ are isomorphic 
where $L_1$ runs through the subroot lattices of all 
reduced presentations. 

The table 5.3 lists for the (conjugacy classes of the) 
exceptions the different ways to write them as 
quasi Coxeter elements of subroot lattices $L_1$, and it lists
the numbers $k_4(L,w)$.

For the other elements, $k_4(L,w)=k_2(L,L_1)$ 
for the unique isomorphism class $(L,L_1)$. 
All these other elements can be found by replacing in the
tables 3.6, 3.2, 3.3 and 3.4 $L_1$ by the possible quasi Coxeter 
elements with subroot lattice of type $L_1$. 
See table 5.1 for the possibilities.
(E.g. $D_5+A_3$ has to be replaced by the two 
possibilities $D_5+A_3$ and $D_5(a_1)+A_3$.)

\medskip
{\bf Table 5.3:}
\begin{eqnarray*}
\begin{array}{l|l|l}
L & \textup{presentation of }w\textup{ as 
quasi Coxeter element} & k_4(L,w)\\ 
 & \textup{in }W(L_1)\textup{ for some subroot lattice }L_1& \\ \hline
G_2 & \www A_2 \sim A_2 & 0\\ \hline
E_7 & D_4(a_1)+2A_1 \sim 2A_3 & 2\\
E_7 & D_4(a_1)+3A_1 \sim 2A_3+A_1 & 1\\
E_7 & D_6(a_1)+A_1 \sim A_7 & 1\\ \hline
E_8 & D_4(a_1)+2A_1 \sim [2A_3]'' & 3\\
E_8 & D_4(a_1)+3A_1 \sim 2A_3+A_1 & 2\\
E_8 & D_5(a_1)+2A_1 \sim D_4+A_3 & 2\\
E_8 & D_6(a_1)+A_1 \sim [A_7]'' & 2\\
E_8 & D_4(a_1)+4A_1 \sim 2A_3+2A_1 & 2\\
E_8 & D_4+D_4(a_1) \sim D_5(a_1)+A_3 & 1\\
E_8 & D_5+A_3 \sim A_7+A_1 \sim D_6(a_1)+2A_1 & 1\\
E_8 & D_6(a_2)+2A_1 \sim 2D_4& 2\\
E_8 & E_6(a_1)+A_2 \sim A_8 & 1\\
E_8 & E_7(a_1)+A_1 \sim D_8 & 1\\
E_8 & E_7(a_3)+A_1 \sim D_8(a_2) & 1 
\end{array}
\end{eqnarray*}

(d) In the case of $F_4$, the following table 5.4 lists 
for (the conjugacy classes of) 
all Weyl group elements all ways to write them
as quasi Coxeter elements of subroot lattices.
See remark \ref{t5.8} (iii) for the notations.
It also lists the numbers $k_4(L,w)$.

\medskip
{\bf Table 5.4:}
\begin{eqnarray*}
\begin{array}{l|l}
\textup{Presentation of }w\textup{ as 
quasi Coxeter element}& k_4(L,w) \\
\textup{in }W(L_1)\textup{ for some subroot lattice }L_1 & \\ \hline
F_4 & 0\\
F_4(a_1) & 0 \\
F_4(a_2)\sim B_4 & 0 \\
F_4(a_3)\sim C_3+A_1 & 0 \\
F_4(a_4)\sim D_4(a_1) \sim B_4(2,2)& 0 \\
B_4(3,1)\sim D_4 & 1 \\
B_4(2,1,1)\sim A_3+\www A_1\sim C_3(2,1)+A_1\sim B_2+2A_1 & 1\\
B_4(1,1,1,1)\sim B_2(1,1)+2A_1\sim 4A_1 & 1 \\
A_2+\www A_2 & 1 \\ \hline
\end{array}
\end{eqnarray*}
\begin{eqnarray*}
\begin{array}{l|l}
w & k_4(L,w) \\ \hline
B_3 & 1 \\
B_3(2,1) \sim A_3 & 1 \\
B_3(1,1,1)\sim 2A_1+\www A_1 & 1 \\
B_2+A_1 \sim C_3(2,1)& 1 \\
B_2(1,1)+A_1 \sim 3A_1 & 2 \\
A_2+\www A_1 & 1 \\
A_1+\www A_2 & 1 \\
C_3 & 1 \\ \hline
\end{array}
\hspace*{0.5cm}
\begin{array}{l|l}
w & k_4(L,w)\\ \hline
A_2 & 2 \\
B_2 & 2 \\
B_2(1,1)\sim 2 A_1 & 2 \\
A_1+\www A_1 & 2 \\
\www A_2 & 2 \\ \hline 
\www A_1 & 3 \\
A_1 & 3 \\ \hline
\emptyset & 4
\end{array}
\end{eqnarray*}
\end{theorem}

{\bf Proof:}
(a) In the cases $A_n,B_n,BC_n,C_n$ and $D_n$, 
any positive cycle in the Weyl group can be written 
as a product $s_{\alpha_1}\circ...\circ s_{\alpha_r}$ where
$\alpha_1,...,\alpha_r$ are roots of the type $\pm e_i\pm e_j$
whose graph $\GG(\{\alpha_1,...,\alpha_r\})$ is a tree.
The subroot lattice $L_1=\sum_{i=1}^r\Z\alpha_i$ is a primitive
sublattice. 

Any element of $W(A_n)$ is a product of positive cycles with disjoint
supports whose union is $\{1,...,n\}$. 
Because the supports are disjoint, the sum of the subroot lattices
of the presentations above of the positive cycles is also a primitive sublattice.
Therefore there \eqref{5.10} holds.

In the cases $B_n$ and $BC_n$, any negative cycle can be written
as a product $s_{\alpha_1}\circ...\circ s_{\alpha_r}\circ s_{\alpha_{r+1}}$
such that $s_{\alpha_1}\circ...\circ s_{\alpha_r}$ is a positive
cycle with graph a tree, and such that $\alpha_{r+1}$ is a short
root which gives a marking of one vertex of the tree.
The subroot lattice $L_1=\sum_{i=1}^{r+1}\Z\alpha_i$ is the primitive
sublattice, which is generated by all the short roots which correspond
to the vertices of the tree.

Any element of $W(B_n)=W(BC_n)$ is a product of positive cycles and/or 
negative cycles with disjoint supports whose union is $\{1,...,n\}$. 
Because the supports are disjoint, the sum of the subroot lattices
of the presentations above of the positive and/or negative 
cycles is also a primitive sublattice. Therefore there \eqref{5.10} holds.

\medskip
(b) In the cases $C_n$ and $D_n$, any pair of negative cycles can 
be written as a product $s_{\alpha_1}\circ...\circ s_{\alpha_a}\circ 
s_{\beta_1}\circ...\circ s_{\beta_b}\circ s_{\alpha_{a+1}}\circ s_{\beta_{b+1}}$
such that $s_{\alpha_1}\circ...\circ s_{\alpha_a}$ and 
$s_{\beta_1}\circ...\circ s_{\beta_b}$ are positive cycles whose graphs 
are disjoint trees and such that $\alpha_{a+1}=e_i-e_j$ and 
$\beta_{b+1}=e_i+e_j$ with $i$ a vertex of one tree and $j$ a vertex of
the other tree. The subroot lattice
$L_1=\sum_{i=1}^{a+1}\Z\alpha_i+\sum_{j=1}^{b+1}\Z\beta_j$ is 
of type $C_{a+b+2}$ respectively $D_{a+b+2}$.

In the case $C_n$, any single negative cycle can 
be written as a product 
$s_{\alpha_1}\circ...\circ s_{\alpha_a}\circ s_{\alpha_{a+1}}$
such that $s_{\alpha_1}\circ...\circ s_{\alpha_a}$ is a positive
cycle with graph a tree, and such that $\alpha_{a+1}$ is a long
root of the type $2e_i$ which gives a marking of one vertex of the tree.
The subroot lattice $L_1=\sum_{i=1}^{a+1}\Z\alpha_i$ is of type
$C_{a+1}$.

Let $w$ be a Weyl group element which is a product of $r$ positive
cycles and $s$ negative cycles with disjoint supports 
whose union is $\{1,...,n\}$. 
One presents the positive cycles as above (in the proof of (a)) 
and as many pairs of negative cycles
as above. At most one (none in the case $D_n$) single negative cycle is left
and is also presented as above.
Let $L_1$ be the subroot lattice of the presentation.
Table 3.1 shows \eqref{5.11}.
One sees easily that no reduced presentation with smaller 
$k_1(L,L_1)$ exists. \eqref{5.12} hold.

\medskip
(c) The tables 7, 9, 10 and 11 in \cite{Ca} list all conjugacy classes
of elements of the Weyl groups of root lattices of the types $G_2,E_6,E_7$
and $E_8$. They give in each case one type of presentation as a
strict quasi Coxeter element. It is easy to find all presentations as
quasi Coxeter elements which are not in the list. One has to find out
which elements in the list are given also by these presentations.
In most cases it is sufficient to compare the characteristic polynomials.
A table of characteristic polynomials is table 3 in \cite{Ca}.

The only cases where this is not sufficient arise for the $E_8$ root lattice
and there for the presentations as quasi Coxeter elements of types
$D_4(a_1)+2A_1$ and $D_6(a_1)+A_1$. In the first case the presentations
as strict quasi Coxeter elements of types $[2A_3]'$ and $[2A_3]''$ have
the same characteristic polynomial, in the second case the presentations
of types $[A_7]'$ and $[A_7]''$.
Because of 
\begin{eqnarray*}
\textup{index}(D_4+2A_1)=2=\textup{index}([2A_3]'')\neq\textup{index}([2A_3]')=1,\\
\textup{index}(D_6+A_1)=2=\textup{index}([A_7]'')\neq\textup{index}([A_7]')=1,
\end{eqnarray*}
lemma \ref{t5.2} (b) tells that $D_4(a_1)+2A_1$ gives the same conjugacy class
as $[2A_3]''$ and that $D_6(a_1)+A_1$ gives the same conjugacy class 
as $[A_7]''$.

\medskip
(d) Table 8 in \cite{Ca} lists 9, 8, 5, 2 and 1 conjugacy classes of
elements of the Weyl group of type $F_4$ of lengths 4, 3, 2, 1 respectively 0.
On the other hand there are 19, 12, 6, 2 and 1 types of presentations
of elements as quasi Coxeter elements of lengths 4, 3, 2, 1 respectively 0:

\begin{eqnarray*}
\begin{array}{l|l}
\textup{length} & \textup{ type of presentation as a quasi Coxeter element}\\ \hline
4 & F_4,F_4(a_1),F_4(a_2),F_4(a_3),F_4(a_4),B_4,B_4(3,1),B_4(2,2),\\
 & B_4(2,1,1),B_4(1,1,1,1), A_3+\www A_1,A_2+\www A_2,C_3+A_1,\\
& C_3(2,1)+A_1,D_4,D_4(a_1), B_2+2A_1,B_2(1,1)+2A_1,4A_1\\ \hline
3 & B_3,B_3(2,1),B_3(1,1,1),B_2+A_1,B_2(1,1)+A_1,A_2+\www A_1,\\
 & A_3,2A_1+\www A_1,A_1+\www A_2, C_3,C_3(2,1),3A_1\\ \hline
2 & B_2,B_2(1,1),\www A_2,A_1+\www A_1,A_2,2A_1\\ \hline
1 & \www A_1,A_1\\ \hline
0 & \emptyset 
\end{array}
\end{eqnarray*}

For those types of presentations as quasi Coxeter elements in the table
above which are not in the table 8 in \cite{Ca}, one has to find
out which conjugacy classes they give. In many cases this is determined
by the characteristic polynomials. The cases where the characteristic
polynomials is not sufficient, can be drawn from lemma 26 in \cite{Ca}.
It lists the presentations as strict quasi Coxeter elements
which give different conjugacy classes, but with the same characteristic
polynomials. Of the 8 pairs in lemma 26 in \cite{Ca}, 
only those 4 are relevant here,
for which presentations as quasi Coxeter elements exist
which are not in table 8 in \cite{Ca}
and which have the same characteristic polynomials. 
These 4 pairs and their characteristic polynomials are as follows:

\begin{eqnarray*}
\begin{array}{l|l|l|l}
D_4 & A_3 & 3A_1 & 2A_1 \\ 
C_3+A_1 & B_2+A_1 & 2A_1+\www A_1 & A_1+\www A_1 \\ \hline
(t^3+1)(t+1) & t^3+t^2+t+1 & (t+1)^3(t-1) & (t+1)^2(t-1)^2
\end{array}
\end{eqnarray*}

The equality
$$s_{e_3}s_{e_4}=s_{e_3+e_4}s_{e_3-e_4}$$
tells that in table 5.4
$$B_2(1,1)\sim 2A_1,\ B_2(1,1)+A_1\sim 3A_1,\ B_3(1,1,1)\sim 2A_1+\www A_1.$$
The equality
$$s_{e_1}s_{\frac{1}{2}(e_1+e_2-e_3-e_4)}s_{\frac{1}{2}(e_1+e_2+e_3+e_4)}
=s_{e_1}s_{e_1+e_2}s_{e_3+e_4}$$
tells that in table 5.4
$$C_3(2,1)\sim B_2+A_1.$$
The equalities
\begin{eqnarray*}
s_{e_2-e_3}s_{e_3}s_{e_4}&=&s_{e_2-e_3}s_{e_3+e_4}s_{e_3-e_4}\qquad
\textup{and}\\
s_{e_1-e_2}s_{e_2-e_3}s_{e_3}s_{e_4}
&=&s_{e_1-e_2}s_{e_2-e_3}s_{e_3+e_4}s_{e_3-e_4}
\end{eqnarray*}
tell that in table 5.4
$$B_3(2,1)\sim A_3\textup{ and }B_4(3,1)\sim D_4.$$
The equivalence $F_4(a_3)\sim C_3+A_1$ in table 5.4 
holds by definition of $F_4(a_3)$.
All other equivalences in table 5.4 follow from comparison of 
characteristic polynomials.\hfill$\Box$

\begin{remarks}\label{t5.11}
(i) From the theorems \ref{t5.10}, \ref{t5.6} and \ref{t3.3} (respectively
the first columns of the tables 3.1--3.6), one can recover
the classification of conjugacy classes of the Weyl group elements
of the root lattices of types $G_2,F_4,E_6,E_7$ and $E_8$ which is
given in the tables 7--11 in \cite{Ca}.

\medskip
(ii) The proof above of theorem \ref{t5.10} had used these tables, but not
in a very crucial way. Those few cases where different conjugacy classes
have the same characteristic polynomials, can be dealt with by hand.
In fact, informations on them are given in the lemmas 26 and 27 in \cite{Ca}.
But theorem \ref{t5.6} on the (strict) quasi Coxeter elements
depends in a crucial way on the results in \cite{Ca}.

\medskip
(iii) The characteristic polynomials of the strict quasi Coxeter elements
in all irreducible root lattices are given in table 3 in \cite{Ca}.

\medskip
(iv) In \cite[(2.3.4)]{Vo} a table similar to table 5.3 
for $E_7$ and $E_8$ is given. But one of the cases for $E_7$ and four
of the cases for $E_8$ are missing there. The case for $E_7$ 
which is missing in \cite[(2.3.4)]{Vo}, is also missing in
\cite[(3.2.9)]{Vo}.
\end{remarks}

\begin{remark}\label{t5.12}
There is a strange correspondence. Define for any 
irreducible root lattice $(L,(.,.),\Phi)$ the two numbers
\begin{eqnarray*}
k_6(L)&:=& |\{\textup{conjugacy classes of quasi Coxeter elements}\}|-1,\\
k_7(L)&:=& |\{\textup{isomorphism classes of pairs }(L,L_1)\textup{ with }L_1\\
&& \textup{ a subroot lattice of full rank with }k_1(L,L_1)=1\}|.
\end{eqnarray*}
Then 
\begin{eqnarray*}
k_6(L)=k_7(L)&&\textup{for }A_n,C_n,D_n,F_4,E_6,E_7,E_8\textup{ and }B_2,\\
&&\textup{but not for }B_n\ (n\geq 3)\textup{ and }G_2,
\end{eqnarray*}
as the following table shows.
\begin{eqnarray*}
\begin{array}{l|l|l|l|l|l|l|l|l|l}
 & A_n & B_n & C_n & D_n & G_2 & F_4 & E_6 & E_7 & E_8 \\ \hline
k_6(L) & 0 & p(n)-1 & \left[\frac{n}{2}\right] & 
\left[\frac{n}{2}\right]-1 & 1 & 4 & 2 & 4 & 8 \\
k_7(L) & 0 & n-1 & \left[\frac{n}{2}\right] & 
\left[\frac{n}{2}\right]-1 & 2 & 4 & 2 & 4 & 8 
\end{array}
\end{eqnarray*}
Here $p(n)$ is the number of partitions of $n$.
\end{remark}

\section{Nonreduced presentations of Weyl group elements}\label{s6}
\setcounter{equation}{0}

\begin{definition}\label{t6.1}
Let $(L,(.,.),\Phi)$ be a p.n. root lattice, and let $w$ be a Weyl
group element. Define the number
\begin{eqnarray}\label{6.1}
k_5(L,w):=\min(k&|& \textup{a presentation }(\alpha_1,...,\alpha_{l(w)+2k})\\
&&\textup{with subroot lattice the full lattice exists}\}. \nonumber
\end{eqnarray}
\end{definition}

Recall the definition \eqref{5.8} of the number $k_4(L,w)$
in the same situation. Let $(\alpha_1,...,\alpha_{l(w)})$ be a reduced
presentation of an element $w$ with subroot lattice $L_1$ such that
$k_2(L,L_1)$ is minimal, i.e. $k_2(L,L_1)=k_4(L,w)=:k$. 
Let $\beta_1,...,\beta_k$
be roots such that $L_1+\sum_{j=1}^k\Z\beta_j=L$. Then obviously
$(\alpha_1,...,\alpha_{l(w)},\beta_1,\beta_1,\beta_2,\beta_2,...,
\beta_k,\beta_k)$ is a presentation with root lattice the full root lattice $L$.
Therefore
\begin{eqnarray}\label{6.2}
k_5(L,w)\leq k_4(L,w).
\end{eqnarray}

\begin{theorem}\label{t6.2}
Let $(L,(.,.),\Phi)$ be a p.n. root lattice, and let $w$ be a Weyl
group element. Then
\begin{eqnarray}\label{6.3}
k_5(L,w)= k_4(L,w).
\end{eqnarray}
\end{theorem}

The proof consists in a reduction to the special case in the following
lemma and in the proof of the following lemma. 
The proof of the lemma is given first.

\begin{lemma}\label{t6.3}
Let $(L,(.,.),\Phi)$ be a p.n. root lattice of some rank $n$, 
and let $w$ be a Weyl group element of length $n-1$. Then
\begin{eqnarray}\label{6.4}
k_5(L,w)=1&\iff& k_4(L,w)=1.
\end{eqnarray}
\end{lemma}

{\bf Proof of lemma \ref{t6.3}:}
If $(L,(.,.),\Phi)$ is reducible with orthogonal summands 
$\bigoplus_{k\in K} L_k$, then $w$ decomposes accordingly
into a product of commuting elements $w_k\in W(L_k)$, and the numbers
$k_4(L,w)$ and $k_5(L,w)$ are additive,
\begin{eqnarray*}
k_4(L,w)=\sum_{k\in K}k_4(L_k,w_k),\qquad k_5(L,w)=\sum_{k\in K}k_5(L_k,w_k).
\end{eqnarray*}
Therefore it is sufficient to prove the lemma and also 
theorem \ref{t6.2}
for the irreducible p.n. root lattices.

Let $(L,(.,.),\Phi)$ be an irreducible p.n. root lattice of rank $n$, and let  
$w$ be a Weyl group element with $l(w)=n-1$.
Then $k_5(L,w)\geq 1$. If $k_4(L,w)=1$ then by \eqref{6.2} also
$k_5(L,w)=1$. Thus it is sufficient to prove $k_5(L,w)=1\Rightarrow k_4(L,w)=1$.

\medskip
The cases ${\bf A_n,D_n,E_6,E_7,E_8}$:
Suppose $k_5(L,w)=1$, and let $(\alpha_1,...,\alpha_{n+1})$ be a presentation
of $w$ whose subroot lattice is the full lattice. 
By theorem \ref{t4.1}, the set $\{\alpha_1,...,\alpha_{n+1}\}$
contains a $\Z$-basis of the full lattice $L$. Using \eqref{2.2b},
we can suppose that $\alpha_1,...,\alpha_n$ is a $\Z$-basis of $L$.
Let $(\beta_1,...,\beta_{n-1})$ be an arbitrary reduced presentation of $w$.
Then
$$s_{\alpha_1}\circ...\circ s_{\alpha_n}=s_{\beta_1}\circ ...\circ s_{\beta_{n-1}}
\circ s_{\alpha_{n+1}}.$$
The subroot lattice of the presentation on the left hand side 
is the full lattice, so it has index one. 
By lemma \ref{t5.2} (b), the index of the subroot lattice of the 
presentation on the right hand side is the same, 
so it is also one. Thus 
$$\sum_{j=1}^{n-1}\Z\beta_j+\Z\alpha_{n+1}=L.$$
This shows here $k_4(L,w)=1.$

\medskip
The cases ${\bf B_n}$ and ${\bf BC_n}$:
Because of theorem \ref{t5.10} (a), $k_4(L,w)=n-l(w)=1$ holds anyway.

\medskip
The cases ${\bf C_n}$:
$k_4(L,w)=1$ holds if and only if a reduced presentation with subroot lattice
of type $A_{n-1}$ or of type $A_{k-1}+C_{n-k}$ for some $k\in\{1,2,...,n-1\}$ 
exists. This follows from table 3.1.
In the case $A_{n-1}$, $w$ is a positive cycle of length $n$.
In the case $A_{k-1}+C_{n-k}$, $w$ is a product of a positive cycle of length $k$
and of one or two negative cycles such that the sum of their lengths is $n-k$.

It rests to show that $w$ is such an element if $k_5(L,w)=1$.
Thus suppose $k_5(L,w)=1$. Let $(\alpha_1,...,\alpha_{n+1})$ be a presentation
of $w$ whose subroot lattice is the full lattice $L$.
By theorem \ref{t4.1}, the set $\{\alpha_1,...,\alpha_{n+1}\}$
contains a $\Z$-basis of the full lattice $L$. Using \eqref{2.2b},
we can suppose that $\alpha_2,...,\alpha_{n+1}$ is a $\Z$-basis of $L$.
Thus $s_{\alpha_2}\circ...\circ s_{\alpha_{n+1}}=:v$ is a quasi Coxeter
element, so either one negative cycle or the product of two negative cycles.

If $\alpha_1$ is a long root, multiplying $v$ from the left with $s_{\alpha_1}$
will turn one of the (one or two) negative cycles into a positive cycle.

If $\alpha_1$ is a short root, so $\alpha_1=\pm e_i\pm e_j$, then the type
of $s_{\alpha_1}\circ v$ depends on the position of the vertices $i$ and $j$
in the supports of the (one or two) negative cycles.
If $i$ and $j$ are in the support of the same negative cycle, then it splits
into two cycles, one positive and one negative.
If $i$ and $j$ are in the supports of different negative cycles, then
$s_{\alpha_1}\circ v$ is a positive cycle of length $n$.

In any case, $w$ is of one of the types which satisfy $k_4(L,w)=1$.

\medskip
The case ${\bf G_2}$: By table 3.6, all subroot lattices of rank 1
are primitive sublattices. Therefore $k_4(L,w)=n-1=1$ holds anyway.

\medskip
The case ${\bf F_4}$: By table 5.4, the only elements $w$ with $l(w)=3$ and 
$k_4(L,w)\geq 2$ are those of type $B_2(1,1)+A_1\sim 3A_1$, 
and the elements of this type satisfy $k_4(L,w)=2$. 
It rests to show for them $k_5(L,w)\geq 2$.

Suppose that such an element $w$ satisfies $k_5(L,w)=1$, and 
let $(\alpha_1,...,\alpha_5)$ be a presentation of $w$ 
whose subroot lattice is the full lattice.
By theorem \ref{t4.1}, the set $\{\alpha_1,...,\alpha_5\}$
contains a $\Z$-basis of the full lattice $L$. Using \eqref{2.2b},
we can suppose that $\alpha_1,...,\alpha_4$ is a $\Z$-basis
of $L$. We may suppose $w=s_{e_1-e_2}s_{e_1+e_2}s_{e_3-e_4}$.
Then
\begin{eqnarray*}
s_{e_1-e_2}s_{e_1+e_2}s_{e_3-e_4}s_{\alpha_5}
=s_{\alpha_1}s_{\alpha_2}s_{\alpha_3}s_{\alpha_4}.
\end{eqnarray*}
Because of the right hand side, this is a quasi Coxeter element
in $W(F_4)$.

First case, $\alpha_5$ is a long root $\alpha_5=\pm e_i\pm e_j$:
Then $\{i,j\}=\{1,2\}$ is impossible because else the
four roots on the left hand side were linearly dependent.
$|\{i,j\}\cap\{1,2\}|=1$ is impossible because
else the element on the left hand side were a Coxeter element
in $W(D_4)$, and this is not a quasi Coxeter element in 
$W(F_4)$. Also $\{i,j\}=\{3,4\}$ is impossible because
else the left hand side were an element of type $4A_1$,
and this is not a quasi Coxeter element in $W(F_4)$, 
or the four roots on the left hand side were linearly dependent.
The first case is impossible.

Second case, $\alpha_5$ is a short root: By conjugation 
and renumbering of the $e_j$ we
can suppose $\alpha_5=\pm e_i$ for some $i$.
Then $i\in\{1,2\}$ is impossible because else the four
roots on the left hand side were linearly dependent.
$i\in\{3,4\}$ is impossible because else the left hand side
were an element of type $B_2+2A_1$, and  this is not a 
quasi Coxeter element in $W(F_4)$. 
The second case is impossible.

Thus $k_5(L,w)\neq 1$, so $k_5(L,w)\geq 2$. This finishes the proof
of the case $F_4$ and the whole proof of lemma 
\ref{t6.3}.\hfill$\Box$

\bigskip
{\bf Proof of theorem \ref{t6.2}:}
Let $(L,(.,.),\Phi)$ be an irreducible p.n. root lattice
of some rank $n$.
At the beginning of the proof of lemma \ref{t6.3} it was
shown that it is sufficient to prove theorem \ref{t6.2} 
in this case. 

Let $w$ be a Weyl group element, and let 
$(\alpha_1,...,\alpha_{l(w)+2k})$ be a presentation with
subroot lattice the full lattice $L$ and with $k=k_5(L,w)$
minimal with this property.

By theorem \ref{t4.1}, the set $\{\alpha_1,...,\alpha_{l(w)+2k})$
contains a $\Z$-basis of the full lattice $L$.
Using \eqref{2.2b}, we can suppose that $\alpha_1,...,\alpha_n$
is a $\Z$-basis of $L$. The element
$s_{\alpha_1}\circ ...\circ s_{\alpha_n}$ has length $n$.
Thus the element
$$v:=s_{\alpha_1}\circ ...\circ s_{\alpha_{n+1}}$$
has length $l(v)=n-1$. And it satisfies $k_5(L,v)=1$.
Lemma \ref{t6.3} applies. Therefore a reduced presentation 
$(\gamma_1,...,\gamma_{n-1})$ of $v$ and a root
$\gamma_{0}$ exist such that $\sum_{i=0}^{n-1}\Z\gamma_i=L$.
Let $L_1$ be the subroot lattice of the presentation
\begin{eqnarray}\label{6.5}
(\gamma_1,...,\gamma_{n-1},\alpha_{n+2},...,\alpha_{l(w)+2k})
\end{eqnarray}
of $w$. As $k$ is minimal, $L_1\subsetneqq L$.
Because of $L_1+\Z\gamma_0=L$, $k_1(L,L_1)=1$.
The presentation of $w$ in \eqref{6.5} shows
$k_5(L_1,w)\leq k-1$. If $k_5(L_1,w)<k-1$ then by adding two times
$\gamma_0$ to a shortest presentation of $w$ with subroot lattice $L_1$,
one obtains also $k_5(L,w)<k$, 
which contradicts the minimality of $k$.
Thus $k_5(L_1,w)=k-1$. Induction on $k$ gives $k_4(L_1,w)=k-1$.
Now
\begin{eqnarray}\label{6.6}
k_4(L,w)\leq k_2(L,L_1)+k_4(L_1,w)=1+(k-1)=k=k_5(L,w).
\end{eqnarray}
Together with \eqref{6.2} this gives \eqref{6.3}.\hfill$\Box$

\section{An application to extended affine root lattices}\label{s7}
\setcounter{equation}{0}

The number $k_5(L,w)$ in definition \ref{t6.1} and theorem \ref{t6.2}
controls existence of quasi Coxeter elements in {\it extended affine
root systems}. These had been defined by K. Saito in \cite[(1.2) and (1.3)]{Sa}.
In \cite{Az} the equivalence with an alternative definition in 
\cite{AABGP} was shown.

The inequalities in lemma \ref{t7.4} below give constraints on a
quasi Coxeter element $w$ in an extended affine root system in terms
of conditions for a nonreduced presentation of the induced element
$\oooo{w}$ in the Weyl group of the associated p.n. root lattice
$L/\Rad(L)$.

\begin{definition}\label{t7.1} 
An {\it extended affine root lattice} is a triple $(L,(.,.),\Phi)$
where $L$ is a lattice, $(.,.):L_\Q\times L_\Q\to \Q$ is a symmetric
positive semidefinite bilinear form, and $\Phi\subset L-\{\alpha\in L\, |\, 
(\alpha,\alpha)=0\}$ is a subset such that the following properties hold.
Here $\langle\beta,\alpha\rangle$ and $s_\alpha$ are defined as in 
\eqref{2.1} and \eqref{2.2}.
\begin{eqnarray}\label{7.1}
\Phi\textup{ is a generating set of }L\textup{ as a }\Z\textup{-module}.\\
\textup{For any }\alpha\in\Phi\ s_\alpha(\Phi)=\Phi.\label{7.2}\\
\langle\beta,\alpha\rangle\in\Z\textup{ for any }\alpha,\beta\in\Z.\label{7.3}
\end{eqnarray}
The elements of $\Phi$ are the {\it roots}, and $\Phi$ is an
{\it extended affine root system}.
\begin{eqnarray}\label{7.4}
W:=\langle s_\alpha\, |\, \alpha\in\Phi\rangle \subset O(L,(.,.))
\end{eqnarray}
is the {\it Weyl group of the extended affine root lattice}.
\end{definition}

\begin{remarks}\label{t7.2}
(i) In \cite{Sa} the definition of an extended affine root system contains 
additionally the following irreducibility property:
\begin{eqnarray}\label{7.5}
\Phi=\Phi_1\cup\Phi_2\textup{ with }\Phi_1\perp\Phi_2
\Rightarrow \Phi_1=\emptyset\textup{ or }\Phi_2=\emptyset.
\end{eqnarray}

\medskip
(ii) Let $(L,(.,.),\Phi)$ be an extended affine root lattice. Because $(.,.)$
is positive semidefinite, the radical of $(L_\R,(.,.))$ is
\begin{eqnarray*}
\Rad(L_\R)
&:=&\{\alpha\in L_\R\, |\, (\alpha,\beta)=0\textup{ for all }\beta\in L_\R\}\\
&=&\{\alpha\in L_\R\, |\, (\alpha,\alpha)=0\}.
\end{eqnarray*}
Define the radicals $\Rad(L):=\Rad(L_\R)\cap L$ 
and $\Rad(L_\Q):=\Rad(L_\R)\cap L_\Q$.
The quotient $L/\Rad(L)$ with the induced bilinear form $(.,.)_{quot}$
and the induced set of roots 
$$\Phi_{quot}:=(\Phi+\Rad(L))/\Rad(L)$$
is obviously a p.n. root lattice. 
It is called the {\it quotient p.n. root lattice}.
Any element $w\in O(L,(.,.))$ induces an element 
$\oooo{w}\in O(L/\Rad(L),(.,.)_{quot})$.
If $w\in W(L)$, then $\oooo{w}\in W(L/\Rad(L))$.
If $\alpha\in\Phi$ induces $\oooo\alpha:=[\alpha]\in L/\Rad(L)$, 
then $\oooo{s_\alpha}=s_{\oooo{\alpha}}\in W(L/\Rad(L))$.

\medskip
(iii) The reducedness property \eqref{2.7} is not required here.
Even if $(L,(.,.),\Phi)$ satisfies it, it
does not necessarily hold for the quotient p.n. root lattice.
That is the reason why in this paper {\it p.n.} root lattices and not only
root lattices are considered.
\end{remarks}

\begin{definition}\label{t7.3}
Let $(L,(.,.),\Phi)$ be an extended affine root lattice of rank $n$.

\medskip
(a) For any element $w$ of its Weyl group, 
a {\it presentation} $(\alpha_1,...,\alpha_k)$,
the {\it length of a presentation}, 
the {\it subroot lattice of a presentation}, and 
the {\it length} $l(w)$ of the element are defined as in definition
\ref{t5.1} (a).

\medskip
(b) An element $w\in W$ is a {\it quasi Coxeter element} if a presentation
of length $n$ exists whose subroot lattice is the full lattice
(this generalizes definition \ref{t5.3} (a)).
\end{definition}

The following simple lemma connects the existence of quasi Coxeter elements
with the numbers $k_5(L/\Rad(L),\oooo{w})$ from section \ref{s6}.
Theorem \ref{t6.2} says $k_5(L/\Rad(L),\oooo{w})=k_5(L/\Rad(L),\oooo{w})$,
and theorem \ref{t5.10} allows to calculate this number.

\begin{lemma}\label{t7.4}
Let $(L,(.,.),\Phi)$ be an extended affine root lattice of rank $n$
with a radical $\Rad(L)$ of rank $r\geq 1$. 
Let $w\in W$ be a quasi Coxeter element, and let $\oooo{w}$ be the induced
element in the Weyl group $W(L/\Rad(L))$ of the quotient p.n. root lattice.
Then
\begin{eqnarray}\label{7.6}
l(\oooo{w})&\leq &n-r\qquad{and}\\
l(\oooo{w})+2k_5(L/\Rad(L),\oooo{w})&\leq& n.\label{7.7}
\end{eqnarray}
\end{lemma}

{\bf Proof:}
\eqref{7.6} is a trivial consequence of \eqref{5.5}, i.e. $l(\oooo{w})=\dim
V_{\neq 1}(\oooo{w})$.
A presentation of length $n$ of $w$ whose subroot lattice is the full lattice $L$
induces a presentation of length $n$ of $\oooo{w}$ whose subroot lattice
is the full lattice $L/\Rad(L)$. This shows \eqref{7.7}.\hfill$\Box$

\begin{examples}\label{t7.5}
(i) The classification of the extended affine root lattices whose
quotient root lattices are inhomogeneous p.n. root lattices is nontrivial,
see \cite{Az} and references therein.

But if $L/\Rad(L)$ is a homogeneous root lattice, then 
there is a sublattice $L_1\subset L$ such that 
$(L_1,(.,.)|_{L_1},\Phi\cap L_1)$ is isomorphic to the quotient
root lattice and $(L,(.,.),\Phi)$ is equal to 
$(L_1\oplus \Rad(L),(.,.),\Phi\cap L_1 +\Rad(L))$.
Thus up to isomorphism, $(L,(.,.),\Phi)$ is determined by the
isomorphism class of the (homogeneous) quotient root lattice 
and by the rank $r$ of the radical.

\medskip
(ii) Let $(L,(.,.),\Phi)$ be an extended affine root lattice of rank $n$
with radical $\Rad(L)$ of rank $r$. Let $w\in W$ be a quasi Coxeter
element such that $\oooo{w}$ has maximal length $l(\oooo{w})=n-r$.
Then lemma \ref{t7.4} and theorem \ref{t6.2} give
\begin{eqnarray}\label{7.8}
r\geq 2k_5(L/\Rad(L),\oooo{w})=2k_4(L/\Rad(L),\oooo{w}).
\end{eqnarray}
\end{examples}

\end{document}